\documentclass{tac}
\usepackage{xy}
\xyoption{all}

\def\mathcaldef#1{\expandafter\def\csname#1\endcsname{{\cal#1}}}

\def\q{\quad}
\def\qq{\quad\quad}
\def\qv{\qq ;\qq}

\def\v{``}

\def\iff{\q\Longleftrightarrow\q}
\def\imp{\Rightarrow}
\def\iimp{\q\Longrightarrow\q}
\def\inc{\hookrightarrow}
\def\sub{\subseteq}
\def\meets{\,\,\cap  \:\!\!\!\! !\,\,\,}
\def\varmeets{\,\,\cap  \negthinspace \negthinspace \negthinspace !\,\,\,}

\newcommand{\tto}{\mathop{\to}\limits}
\newcommand{\ttto}{\mathop{\longrightarrow}\limits}
\newcommand{\ppi}{\mathop{\pi}\limits}

\def\df{\overleftarrow}
\def\dof{\overrightarrow}
\def\dbf{\overline}
\def\ul{\underline}

\def\down{\downarrow\!\!}
\def\up{\uparrow\!\!}
\def\rup{\!\!\uparrow}
\def\rdown{\!\!\downarrow}

\def\la{\langle}
\def\ra{\rangle}
\def\adj{\dashv}
\def\op{^{\rm op}}
\def\inv{^{-1}}
\def\d{{\rm D}}
\def\n{{\rm N}}
\def\z{{\rm Z}}
\def\ee{{\rm E}}
\def\pb{{\rm P}}
\def\CatX{\Cat\!/\! X}
\def\CatY{\Cat\!/Y}
\def\CatFX{\Cat\!/\!{\cal F}\!X}
\def\GrphX{\Grph\!/\! X}
\def\comp{\Gamma_{!}}
\def\disc{\Gamma^{*}}
\def\pt{\Gamma_{*}}
\def\eps{\varepsilon}

\mathcaldef{P}
\mathcaldef{A}
\mathcaldef{B}
\mathcaldef{C}
\mathcaldef{D}
\mathcaldef{N}
\mathcaldef{Z}
\mathcaldef{Y}

\mathbfdef{Set}
\mathbfdef{Cat}
\mathbfdef{Pos}
\mathbfdef{Grph}

\mathbfdef{Rel}
\mathbfdef{Eq}
\mathbfdef{1}
\mathbfdef{2}
\mathbfdef{3}
\def\Lendo{l{\bf Eset}}
\def\Rendo{r{\bf Eset}}
\def\Lidem{l{\bf Idm}}
\def\Ridem{r{\bf Idm}}
\mathbfdef{Bij}
\mathbfdef{y}
\mathbfdef{e}
\mathbfdef{i}
\mathbfdef{r}

\mathrmdef{hom}
\mathrmdef{ten}
\mathrmdef{true}
\mathrmdef{false}
\mathrmdef{ev}
\mathrmdef{id}
\mathrmdef{and}
\def\andq{\q\and\q}
\mathrmdef{Nat}
\mathrmdef{Dinat}
\mathrmdef{Lim}
\mathrmdef{Colim}
\mathrmdef{Cone}
\mathrmdef{Wedge}
\mathrmdef{End}
\mathrmdef{Coend}
\mathrmdef{Cone}
\mathrmdef{gcd}
\mathrmdef{lcm}
\mathrmdef{fix}
\def\no{{\rm not}}

\let\pf\proof
\let\epf\endproof
\def\eq{\begin{equation}}
\def\eeq{\end{equation}}
\def\eqa{\begin{eqnarray}}
\def\eeqa{\end{eqnarray}}
\def\eqa*{\begin{eqnarray*}}
\def\eeqa*{\end{eqnarray*}}
\def\nn{\nonumber}

\newtheorem{prop}{Proposition}
\newtheorem{lemma}{Lemma}
\newtheorem{corol}{Corollary}

\author{Claudio Pisani}
\address{via Gioberti 86,\\ 10128 Torino, Italy.}
\title{Components, complements \\ and reflection formulas}
\copyrightyear{2006}
\keywords{categories over a base, discrete fibrations and opfibrations, reflections and coreflections, 
components, tensor, negation and complements, atoms, strong dinaturality, limits and colimits, ends and coends,
Cauchy completion, Kan extensions, graphs and evolutive sets}
\amsclass{18Axx}
\eaddress{pisclau@yahoo.it}

\begin{document}

\maketitle

\begin{abstract}
Some basic features of the simultaneous inclusion of discrete fibrations and discrete opfibrations 
in categories over a base category $X$ are considered. 
In particular, we illustrate the formulas 
\[ (\down P)x = \ten(x/X,P) \qv (P\rdown)x = \hom(X/x,P) \]
which give the reflection $\down P$ and the coreflection $P\rdown$ of a category $P$ over $X$ in discrete fibrations. 
The explicit use of the \v tensor functor" 
\[ \ten := \comp(-\times -):\CatX\times\CatX \to \Set\]
given by the components of products, allows a vast generalization of the corresponding analysis in the two-valued context.
For any df $A$, the functor $\ten(A,-):\CatX \to \Set$ has a right adjoint $\neg A$ valued in dof's (and vice versa);
such a complement operator, which in the two-valued case reduces to the classical complementation between lower and upper parts 
of a poset, turns out to be an effective tool in the set-valued context as well.   


Various applications of the formulas and of the accompanying conceptual frame are presented.
\end{abstract}
\tableofcontents

\newpage
{\footnotesize
\begin{enumerate}
\item
Introduction 
\item
The two-valued context 
\begin{enumerate}
\item
Truth values, posets, lower and upper parts. (\ref{2sub1})
\item
The formulas via Kan extensions. (\ref{2sub2}) 
\item
The \v meets" predicate and complements. (\ref{2sub3})
\item
Atoms in a poset. (\ref{2sub4})
\item
Negation and complements in $\P X$. (\ref{2sub5})
\item
The reflection formula via duality. (\ref{2sub6})
\item
Properties of \v meets" and of complements in $\P X$. (\ref{2sub7})
\item
The reflection formula again. (\ref{2sub8})
\end{enumerate}
\item
Categories over a base and components 
\begin{enumerate}
\item
Categories over a base and discrete fibrations. (\ref{3sub1})
\item
Strong dinaturality. (\ref{3sub2})
\item
Components of graphs and categories. (\ref{3sub3})
\item
Limits and colimits of set functors. (\ref{3sub4})
\item
Ends and coends. (\ref{3sub5})
\end{enumerate}
\item
The tensor functor and complements of categories over a base 
\begin{enumerate}
\item
The tensor functor and complements. (\ref{4sub1})
\item
Atoms. (\ref{4sub2})
\item
The tensor functor on $\CatX$ and tensor products. (\ref{4sub3})
\item
Negation and complements in $\CatX$. (\ref{4sub4})
\item
Properties of the tensor functor and of complements in $\CatX$. (\ref{4sub5})
\end{enumerate}
\item
The reflection and the coreflection of categories over a base in discrete fibrations  
\begin{enumerate}
\item
The reflections of objects and the Yoneda Lemma. (\ref{5sub1})
\item
The coreflection and reflection formulas. (\ref{5sub2})
\item
The formulas work. (\ref{5sub3})
\end{enumerate}
\item
Applications 
\begin{enumerate}
\item
Kan extensions via reflection. (\ref{6sub1})
\item
Limits and colimits. (\ref{6sub2})
\item
Atoms in $\CatX$ and their reflections. (\ref{6sub3})
\item
Graphs and evolutive sets. (\ref{6sub4})
\end{enumerate}
\end{enumerate}
}

\section{Introduction}
\label{introduction}

We introduce the subject by a simple example.
Consider an equivalence relation $\sim$ on a set $X$. 
(To fix ideas, let $X$ be the plane and let the equivalence classes of $\sim$ be given by a tiling of $X$.)
Then, among the parts of $X$ we have the closed parts, that is those subsets $P\sub X$ such that
\[ x\in P \andq x\sim y \iimp y\in P \]
so that $P$ is an union of equivalence classes.
The inclusion $i$ of closed parts in all the parts of $X$, has both left and right adjoints
\[ \dbf{(-)} \adj i \adj \ul{(-)} \]
giving the \v best" outer and inner approximations of a part $P$ by closed parts:
\[ \ul P \sub P \sub \dbf P \]
Note that for any singleton part $x\in X$, $\dbf x$ is the equivalence class of $x$.
It is intuitively clear that the coreflection and the reflection are given explicitly by
\eq \label{eq1}
x \in \ul P \iff \dbf x \sub P 
\eeq
\eq \label{eq2}
x \in \dbf P \iff \dbf x \meets P 
\eeq
where $\meets$ is the \v meets" predicate: $P$ meets $Q$ if they have a non-void intersection:
\[ P \meets Q \iff P \cap Q \not\sub \emptyset \]
(warning: the similarity of this symbol with the cotensor operator $\,\cap \:\!\!\!\! |\,$ is fortuitous). 

Dropping the symmetry condition, let $X$ be a poset. Now, a part $P$ of $X$ may be donward or upward closed.
So we have the inclusions
\[ i : \df X \inc \P X \qv j : \dof X \inc \P X \]
of lower and upper parts in all the parts of $X$, which have left and right adjoints
\[ \down(-) \adj i \adj (-)\rdown \qv \up(-) \adj j \adj (-)\rup \]
explicitly given by
\eq \label{eq3}
x \in P\rdown \iff \down x \sub P \qv x \in P\rup \iff \up x \sub P \\
\eeq
\eq \label{eq4}
x \in\,\down P \iff \up x \meets P \qv x \in\,\up P \iff \down x \meets P 
\eeq
which are the non-symmetrical version of the~(\ref{eq1}) and~(\ref{eq2}) above. 

In Section~\ref{posets} we prove the formulas~(\ref{eq3}) and~(\ref{eq4}) in such a way that they can be generalized
in a natural way to the set-valued context. 
In fact, after preparing the setting in Section~\ref{cats} and Section~\ref{complements}, we show in Section~\ref{proof} 
that the coreflection and the reflection of a category $P$ over $X$ in df's and dof's over $X$ have the following form 
\eq \label{eq5}
(P\rdown)x = \hom(\down x,P) \qv (P\rup)x = \hom(\up x,P) \\
\eeq
\eq \label{eq6}
(\down P)x = \ten(\up x,P) \qv (\up P)x = \ten(\down x,P)
\eeq
which are exactly the set-valued versions of~(\ref{eq3}) and~(\ref{eq4}).
(Note that $\down x \cong X/x$ and $\up x \cong x/X$ correspond to the representable presheaves.)
While the former are immediate, the latter rest on the (partially defined) \v complement operator" given by 
exponentials on discrete basis:
\[ \neg:=(\disc -)^- : (\CatX)\op\times\Set \to \CatX \]
which takes a df $A$ in the functor $\neg A : \Set \to \CatX$ valued in dof (and vice versa) 
and allows to prove the almost obvious two-valued formula
\eq
\down P \meets Q \iff P \meets\!\up Q
\eeq
in the set-valued context too:
\eq
\ten(\down P,Q) \cong \ten(P,\up Q)
\eeq
That the formulas~(\ref{eq5}) and~(\ref{eq6}) actually give the desired coreflection and reflection is also 
checked in Section~\ref{proof}.

Among the various applications of the theory discussed in Section~\ref{applications} there are:
\begin{itemize}
\item
A treatment of limits and colimits which makes clear the relevance of the reflection in df's.
\item
An analysis of (weak) atoms in $\CatX$, that is those categories $x$ over $X$ such that
\eq
\hom(x,A) \cong \ten(x,A)
\eeq
for any df or dof $A$. Among them we find not only the objects of $X$, but also the idempotent arrows in $X$.
The reflections of atoms are, as presheaves on $X$, the retracts of representable functors; 
so they generate the Cauchy completion of $X$, displaying it as the Karoubi envelope of $X$.
\item
A proof of the Kan extensions formulas for set functors.
\item
Applications of the formulas in \v other" contexts; e.g. we show how they give in a very pictorial way the reflection
of graphs in evolutive sets. 
In fact, one of the themes of the present work is that of treating on an equal footing the coreflections and the reflections: 
while the former are easily analyzed by figures via the $\hom$ functor, 
a similar direct analysis is now possible for the latter using the tensor functor instead.
The symmetry or duality of the two approaches is due to the fact that we use atomic shapes.  
\end{itemize}
The present paper has evolved from the preprint~\cite{bip} and the presentation of a part of it at the 
International Conference on Category Theory (CT06) held at White Point in June 2006.
While~\cite{bip} aimed toward a more general setting, we here develop several new items.
Among them there is a natural presentation of strong dinaturality and a corresponding
analysis of ends and coends in Section~\ref{cats}; 
a more careful study of the \v duality" in $\CatX$ in Section~\ref{complements},
which is obtained by isolating some key features of the classical duality in $\P X$ 
and between sieves and cosieves of a poset in Section~\ref{posets};
and the treatment of limits in Section~\ref{applications}, which gives a clearer perspective on some points
studied in~\cite{par73}.

It is the author's opinion that basic category theory still needs clarification, both for divulgation purposes and
as a basis for further research. We hope to show that a careful analysis of the dialectical relationships between categories over
a base and discrete fibrations, along with an explicit consideration of the tensor functor, is a step in this direction.  

\noindent {\em Acknoweledgements.} Useful discussions with R.J. Wood and Vincent Schmitt are gratefully acknoweledged.
I also acknoweledge the deep influence of the work of Lawvere and of Par\'e (in particular~\cite{par73}). 
My general attitude in preparing the present work was that of \v noticing, cultivating, and teaching helpful examples
of an elementary nature" in category theory, as advocated in~\cite{law86}.

\section{The two-valued context}
\label{posets}

In this section we derive in several ways the formulas for the reflection and coreflection of 
the parts of a poset in lower and upper parts.
In so doing, we introduce the conceptual frame that will serve for the corresponding analysis
in the set-valued context as well: the reflection and coreflection of categories over a base
in discrete fibrations. 

Note though that we are not claiming that there is a general enriched theory subsuming the
two-valued and set-valued cases as particular instances. On the contrary, the question of
the existence of such a theory is one of the points that may deserve further research.

\subsection{Truth values, posets, lower and upper parts}

\label{2sub1}
In this section we consider $\2$ as the truth values category, and so denote its arrows by
\[ \false \vdash \false \qv \false \vdash \true \qv \true \vdash \true \]
The category $\2$ is cartesian closed, the product being \v and", and exponentials being all
equal to $\true$ except
\[ \true \imp \false\q = \q\false \]
Following~\cite{law73}, a poset $X$ is a $\2$-category, that is 
\[ x\leq y \iff X(x,y) \in \2 \]
is the truth value of the dominance of $x$ by $y$.
(Note that $\iff$ is a way to express equality in $\2$.) 

A functor $A:X\op\to\2$ corresponds to a subset $A\sub |X|$ of the underlying set of $X$ by
\[ x\in A \iff Ax \]
Such a subet is a \v lower set" or \v sieve" or \v order ideal":  
\[ x\in A \andq y\leq x \q\vdash\q y\in A \]
Dually, a functor $A:X\to\2$ corresponds to an \v upper set" or \v cosieve": 
\[ x\in A \andq x\leq y \q\vdash\q y\in A \]
If we denote by 
\[ \df X = \2^{X\op} \qv \dof X = \2^X \]
the categories of lower and upper sets, we have the inclusions
\[ i : \df X \inc \P X \qv j : \dof X \inc \P X \]
that have left and right adjoints:
\[ \down(-) \adj i \adj (-)\rdown \qv \up(-) \adj j \adj (-)\rup \]
Observe that each singleton part $x\in X$ has reflections 
given by the representable (or principal) sieve $\down x = X(-,x)$
\[ y\in\,\down x\iff y\leq x\]
and the representable cosieve $\up x = X(x,-)$ 
\[ y\in\,\up x\iff x\leq y\]
that is (Yoneda lemma)
\[ x\in A \iff x\sub A \iff \down x\sub A \qv x\in D \iff x\sub D \iff \up x\sub D \]
for any lower part $A$ and any upper part $D$.

\subsection{The formulas via Kan extensions}
\label{2sub2}

If we denote by 
\[ \disc \adj \pt : \Cat\to\Set \]
the \v discrete" and \v objects" functors (see section~\ref{3sub3}), the counit of the adjunction 
\[ |X| = \disc\pt X \to X \] 
is the functor that includes the discrete category associated to $X$ in $X$ itself. 
Note that $|X\op| = |X|$. In our case, we can pose
\[ \P X = \2^{|X|} = \2^{|X\op|} = \P X\op \]
and consider the above inclusions as induced by $|X|\to X$ and $|X|\to X\op$,
so that the left and right Kan extensions along these functors give their left and right adjoints.

Thus, the reflection $\down(-)$ can be obtained by the coend formula for left Kan extensions:
\[ (\down P)x = \Coend_{|X|}(X(x,-)\,\,\and\,\, P) \qv (\up P)x = \Coend_{|X|}(X(-,x)\,\,\and\,\, P) \]
that is 
\[ x\in\,\down P \iff \exists y\, (x\leq y \,\,\and\,\, y\in P) \qv x\in\,\up P \iff \exists y\, (y\leq x \,\,\and\,\, y\in P) \]
So the coend in question (the tensor product of functors on $|X|$) reduces to the \v meets" predicate $\varmeets$ on $\P X$
($P$ meets $Q$ iff they have a non-void intersection) giving the reflection formula~(\ref{eq4}):
\[ x \in\,\down P \iff \up x \meets P \qv x \in\,\up P \iff \down x \meets P \]
Similarly the coreflection is given by the end formula:
\[ (P\rdown)x = \End_{|X|}(X(-,x)\imp P) \qv (P\rup)x = \End_{|X|}(X(x,-)\imp P) \]
that is
\[ x\in P\rdown \iff \forall y\, (y\leq x \,\imp\, y\in P) \qv x\in P\rup \iff \forall y\, (x\leq y \,\imp\ y\in P) \]
so that the end in question reduces to the inclusion predicate on $\P X$, giving the coreflection formula~(\ref{eq3}): 
\[ x \in P\rdown \iff \down x \sub P \qv x \in P\rup \iff \up x \sub P \]
This method would fail in the set-valued context, since our main concern, the inclusion of df's in categories over a base
is not a Kan-like functor (note e.g. that $\CatX$ is not a presheaf category).

\subsection{The coreflection formula via figures}
\label{2sub3}

One can obtain the coreflection formula in a more straightforward manner by using the singleton parts as figures:
\[
\begin{array}{c}
x\in P\rdown \\ \hline
x\sub P\rdown \\ \hline 
\down x\sub P\rdown \\ \hline 
\down x\sub P 
\end{array}
\qq\qq
\begin{array}{c}
x\in P\rup \\ \hline
x\sub P\rup \\ \hline 
\down x\sub P\rup \\ \hline 
\down x\sub P 
\end{array}
\]
where only the adjunction laws have been used.

\subsection{The \v meets" predicate and complements}
\label{2sub4}

To derive the reflection formula it seems natural to use duality.
This can be done in a simple-minded way (see section~\ref{2sub7}) or in a more \v wise" manner that allows a similar 
derivation in the set-valued context (see section~\ref{2sub9}). 
Duality has two aspects: complementation, that is related to exponentiation, and negation, 
that is complementation in the truth-values category which lifts to functors valued therein. 
Of course the two aspects are strictly related.

Let $\A$ be any bounded poset, that is a poset with a maximum (a terminal object) $\top$ and a minimum (an initial object) $\bot$.
There are functors
\[ \comp \adj \disc \adj \pt :\A\to\2 \] 
given by
\begin{enumerate}
\item
$\disc \true = \top \qv \disc \false = \bot$
\item
$\pt x = \true {\q\rm iff\q} \top\leq x$
\item
$\comp x = \false {\q\rm iff\q} x\leq\bot$
\end{enumerate}
So we may interpret $\pt$ as giving the truth value of the proposition ``$x$ is a maximum of $\A$",
and $\comp$ as giving the truth value of the proposition ``$x$ is not a minimum of $\A$" or ``$x$ is not empty".

If $\A$ is also a meet semilattice (that is, it has products) then it is natural to define the \v meets" predicate as follows: 
\[ x\meets y := \comp (x\wedge y) \] 
that is \v $x$ meets $y$ iff their meet is not empty". 
\begin{definition}
If the functor $x\varmeets - :\A\to\2$ has a right adjoint, we call it the \v complement" of $x$ and denote it by $\neg x$: 
\[ x\meets -\adj\neg x :\2\to\A \] 
\end{definition}
To motivate the above definition, recall that if the exponential $x\imp y$, characterized by 
\( \q z\leq x\imp y \iff z\wedge x \leq y \,\), 
exists then it is called the \v relative pseudocomplement" of $x$ with respect to $y$. 
If $\,\neg x := x\imp\bot$ exists in the bounded semilattice $\A$ then it is called the \v pseudocomplement" of $x$. 
\begin{prop} \label{24prop1}
An element $x$ of a bounded semilattice $\A$ has complement $\neg x:\2\to\A$ iff the pseudocomplement $x\imp\bot$ exists.
If this is the case, $(\neg x)\false$ is the pseudocompement of $x$ itself, while $(\neg x)\true$ 
is a obviously a maximum of $\A$.
\end{prop}
\pf
In one direction, if $x\imp\bot$ exists (and since $x\imp\top = \top$ always exists), the desired right adjoint 
can be obtained by composition of adjoints: $\neg x = x\imp\disc -$. In the other direction 
\[
\begin{array}{c}
y\leq (\neg x)\false \\ \hline
x\meets y \vdash \false \\ \hline 
x\wedge y \leq \bot 
\end{array}
\]
which gives the required adjunction law.
\epf
So the pseudocomplement $\neg x$ is only a trace of a functor $\neg x : \2\to\A$, namely its value at $\false\in\2$. 
The fact that the other value is $(\neg x)\true = \top$ conceals its significance, making it superfluous in the two-valued context.
On the contrary, it becomes essential in the set-valued context.

\subsection{Atoms in a poset}
\label{2sub5}

Let $\A$ be a bounded semilattice. The \v meets" predicate therein allows a very natural definition of atom:
an element $x\in\A$ is an \v atom" of $\A$ when for any $y\in\A$
\[ x\meets y \iff x\leq y \]
That is, if we see the elements of $\A$ as \v parts", an atom is \v so small" to be included in any part that it meets, 
but also \v big enough" to meet any part in which it is included (in particular it meets itself).
The second condition excludes the bottom $\bot$.
\begin{example}
Let $\A$ be the poset of positive natural numbers ordered by divisibility. 
Since $\bot = 1$ and $n\wedge m =\gcd(n,m)$, two numbers meet each other in $\A$ iff they have a common divisor greater than $1$.
Thus the atoms of $\A$ are the prime numbers. 
So in this case the atoms coincide with the prime elements, that is those elemets such that the representable $\q x\leq -\q$
preserves the finite sups (i.e. the finite \lcm's).  
\end{example}
\begin{prop}  \label{25prop1}
If an atom of $\A$ has complement then it is a completely prime element.
In particular, the atoms of a Heyting algebra are completely prime.
\end{prop}
\pf
Indeed, the complement $\neg x:\2\to\A$ of the atom $x\in\A$ is the right adjoint of $\q x\meets -\q =\q x\leq -\q$
so that $\q x\leq -\q$ preserves all sups (colimits).
\epf
\begin{prop}  \label{25prop2}
The singleton parts $x\in X$ are atoms in $\P X$:
\[ x\meets P \iff x\sub P \iff x\in P \]
for any part $P\in\P X$.
\end{prop}
\epf

\subsection{Negation and complements in $\P X$}
\label{2sub6}

If $\A$ is a bounded meet semilattice, pseudocomplementation
\[ -\imp\bot:\A\to\A\op \]
is in general only a partially defined functor.
If $\A=\2$, pseudocomplementation is usually called \v negation"
\[ \no := -\imp\false:\2\to\2\op \]
and is an isomorphism which induces
\[ \no^{X\op} : \2^{X\op}\to(\2\op)^{X\op} \cong (\2^X)\op \qv \no^X : \2^X\to(\2\op)^X \cong {(\2^{X\op})}\op \] 
giving isomorphisms
\[ \neg : \df X\ttto(\dof X)\op \qv \neg : \dof X\ttto(\df X)\op \] 
Explicitly, for any $A:X\op\to\2$, $\neg A$ is (the opposite of) $\no\circ A$, and similarly for any $D:X\to\2$: 
\eq \label{26eq1}
\neg A : x \mapsto \no\, Ax  \qv  \neg D : x \mapsto \no\, Dx
\eeq
In $\P X$, relative pseudocomplemets (exponentials) are computed \v pointwise":
\[ (P\imp Q)x \iff Px\imp Qx \]
In particular, pseudocomplementation in $\P X$
\eq \label{26eq2}
(-)' := -\imp\emptyset : \P X \ttto (\P X)\op 
\eeq 
is the isomorphism given by
\eq \label{26eq3}
P' : x \,\,\mapsto\q Px\imp\false \q = \,\,\no\, Px 
\eeq
By comparing~(\ref{26eq1}) and~(\ref{26eq3}) we have the following commutative diagrams, whose rows are isomorphisms:
\eq \label{26eq4}
\xymatrix{
{\df X} \ar[rr]^-{\neg}\ar@{_{(}->}[dd]^i     && {(\dof X)\op}\ar@{_{(}->}[dd]^{j\op} \\ \\
{\P X}      \ar[rr]^-{(-)'}                   && ({\P X})\op }
\qq\qq
\xymatrix{
{\dof X} \ar[rr]^-{\neg}\ar@{_{(}->}[dd]^j && {(\df X)\op}\ar@{_{(}->}[dd]^{i\op} \\ \\
{\P X}      \ar[rr]^-{(-)'}                && ({\P X})\op }
\eeq
Each one of the diagrams~(\ref{26eq4}) can be obtained from the other by substituting $X$ with $X\op$ 
(the lower row remains unchanged, since $\P X = \P(X\op)$);
furthermore each one can be merged with the opposite of the other giving
\eq \label{26eq5}
\xymatrix{
{\df X} \ar[rr]^-{\neg}\ar@{_{(}->}[dd]^i  && {(\dof X)\op}\ar@{_{(}->}[dd]^{j\op} \ar[rr]^-{\neg\op}       && {\df X}\ar@{_{(}->}[dd]^i \\ \\
{\P X}      \ar[rr]^-{(-)'}                && ({\P X})\op                          \ar[rr]^-{(-)'^{\rm op}} && {\P X}} 
\qq
\xymatrix{
{\dof X} \ar[rr]^-{\neg}\ar@{_{(}->}[dd]^j  && {(\df X)\op}\ar@{_{(}->}[dd]^{i\op} \ar[rr]^-{\neg\op}       && {\dof X}\ar@{_{(}->}[dd]^j \\ \\
{\P X}      \ar[rr]^-{(-)'}                 && ({\P X})\op                         \ar[rr]^-{(-)'^{\rm op}} && {\P X}} 
\eeq
whose horizontal compositions are identities.
For future reference, we summarize the above considerations in the following proposition
\begin{prop} \label{26prop1}
The pseudocomplement of a lower part is an upper part, and conversely. Equivalently, the complement of a lower part 
is valued in upper parts, and conversely. 
\end{prop}
\epf
\begin{remark}
More generally, note that if $A$ and $D$ are a lower and an upper part, then their exponential 
$A\imp D$ in $\P X$ is an upper part (and dually $D\imp A$ is a lower part).
Indeed, since $\P X$ is in fact a boolean algebra, $A\imp D$ is the union of two upper parts
\[ A\imp D = A'\cup D \]
Alternatively, observe that, since exponentials are computed \v pointwise", if $x\leq y$
\[ (A\imp D)x \iff Ax\imp Dx \iimp Ay\imp Dy \iff (A\imp D)y \]
Remarkably, this holds also for discrete fibrations and discrete opfibrations, when the exponential
is taken in $\CatX$ (see Proposition~\ref{31prop1}).
\end{remark}

\subsection{The reflection formula via duality}
\label{2sub7}

The duality summarized by diagrams~(\ref{26eq4}) and~(\ref{26eq5}) can be used to reduce the reflection to the 
known coreflection of section~\ref{2sub3}.
Indeed, in general a left adjoint of a functor $F:\A\to\B$ is a right adjoint of $F\op:\A\op\to\B\op$.
So to find the left adjoint to $j : \dof X \inc\P X$ we look for the right adjoint to $j\op : (\dof X)\op \inc(\P X)\op$.
But by diagram~(\ref{26eq4}) the latter is the \v same" as $i : \df X \inc \P X$ via the isomorphisms $\neg$ and $(-)'$.
Explicitly, we have $\up P = \neg(P'\rdown)$, that is
\[
\begin{array}{c}
x\in\,\up P \\ \hline
x\in \neg(P'\rdown) \\ \hline 
\no\,\,\, x\in P'\rdown \\ \hline 
\no\, \down x\sub P' \\ \hline 
\no\,\,\, \down x \cap P \sub \emptyset \\ \hline 
\down x\meets P
\end{array}
\]
This method too would fail in the set-valued context, since we have no such duality between df's and dof's.

\subsection{Properties of \v meets" and complements in $\P X$}
\label{2sub8}

We now prove an \v adjunction-like" property of the meets predicate on $\P X$ 
which is based on the complement functor and allows us to derive directly the reflection formula.
\begin{prop} \label{28prop1}
Let $X$ be a poset, $A$ a lower part and $D$ an upper part of $X$. Then for any part of $P$ of $X$
\[  P\meets A  \iff \up P\meets A \qv P\meets D \iff \down P\meets D \]
\end{prop}
\pf
We prove the first one, the other being symmetrical:
\[
\begin{array}{c}
P\meets A \vdash \false \\ \hline 
P \sub (\neg A)\false \\ \hline 
\up P \sub (\neg A)\false \\ \hline 
\up P\meets A \vdash \false
\end{array}
\]
where we have used Proposition~\ref{24prop1} and Proposition~\ref{26prop1}.
\epf

The following proposition is almost obvious, but we record it since it will be generalized in the set-valued context:
\begin{prop} \label{28prop2}
If $A$ and $B$ are both lower sets (or both upper sets), then 
\[ A\leq B \iff \neg B \leq \neg A \]
where the right-hand condition may be seen in $\P X$ or equivalently in $(\P X)^\2$, as a natural transformation between 
the complements $\neg B, \neg A:\2\to\P X$.
Briefly, the \v strong contraposition law" holds.
Furthermore, lower parts are properly analyzed by upper parts (and conversely) via the \v meets" predicate:
\[ A\leq B \iff A\meets - \,\leq\, B\meets- \]
where the right-hand side refers to functors $\dof X\to\2$: 
\[ \forall D\in\dof X \quad (A\meets D\q\imp\q B\meets D) \]
\end{prop}
\epf

\subsection{The reflection formula again}
\label{2sub9}

We are now in a position to derive the reflection formula in a way straightforwardly generalizable to the set-valued context.
Since singleton parts are atoms in $\P X$ (see Proposition~\ref{25prop2}), we can analyze parts $P$ 
not only with figures like in section~\ref{2sub3} 
but also with the \v meets" predicate:
\[
\begin{array}{c}
x\in\,\,\down P \\ \hline
x\meets\!\down P \\ \hline 
\up x\meets\!\down P \\ \hline 
\up x\meets P 
\end{array}
\qq\qq
\begin{array}{c}
x\in\,\,\up P \\ \hline
x\meets\!\up P \\ \hline 
\down x\meets\!\up P \\ \hline 
\down x\meets P 
\end{array}
\]
where we used twice the adjunction-like properties of Proposition~(\ref{28prop1}). 
We stress the similarity between the derivations of the reflection and coreflection formulas by posing them side by side:
\[
\begin{array}{c}
x\in\,\,\down P \\ \hline
x\meets\!\down P \\ \hline 
\up x\meets\!\down P \\ \hline 
\up x\meets P 
\end{array}
\qq\qq
\begin{array}{c}
x\in P\rdown \\ \hline
x\sub P\rdown \\ \hline 
\down x\sub P\rdown \\ \hline 
\down x\sub P 
\end{array}
\]

\section{Categories over a base and components}
\label{cats}

We collect some facts that will be needed in the sequel. But in so doing we hopefully also cast some new light 
on categorical concepts such as strong dinaturality, components, coends and exponentials in $\CatX$.

\subsection{Categories over a base and discrete fibrations}
\label{3sub1}

We denote by $\df X$ and $\dof X$ the full subcategories of $\CatX$ of discrete fibrations and discrete opfibrations:
\eq
\xymatrix@1@C=3pc{{\df X}\,\, \ar@{^{(}->}[r]^-i & \CatX & \,{\dof X} \ar@{_{(}->}[l]_-j }
\eeq
By a well-known argument, we have equivalences
\eq  \label{31eq1}
\Set^{X\op} \simeq \df X \qv \Set^X \simeq \dof X 
\eeq
and so also full and faithful functors
\eq  \label{31eq2}
\xymatrix@1@C=3pc{ \Set^{X\op}\, \ar@{^{(}->}[r] & \CatX & \,\,\Set^X \ar@{_{(}->}[l] }
\eeq 
When opportune, we do not distinguish notationally a df and the corresponding presheaf; e.g., if $A$ is a presheaf 
we may say \v $A$ over $X$" to emphasize that we are thinking of it as a category over $X$ 
via the above equivalence (the \v category of elements").

Any object $x:\1\to X$ of the base category $X$ is a category over $X$. 
If $\pi:P\to X$ is a category over $X$, the product $P\times x$ in $\CatX$ has as total category the fiber $Px$ over $x$.
Similarly, for any arrow $f:\2\to X$ in $X$, we have a \v fiber" $Pf$ which is in fact a category over $\2$,
as displayed by the right-hand pullback:
\eq 
\xymatrix@R=3pc@C=3pc{
Px   \ar[r]\ar[d]    & P \ar[d]^{\pi} \\
{\1} \ar[r]^x        & X }
\qq\qq
\xymatrix@R=3pc@C=3pc{
Pf   \ar[r]\ar[d]    & P \ar[d]^{\pi} \\
{\2} \ar[r]^f        & X }
\eeq
Exponentials in $\CatX$ may not exist (for example, $x:\1\to X$ is not exponentiable if there is another object $y\cong x$ in $X$;
see e.g.~\cite{bun00} for an account on exponentiability in $\CatX$). 
When they do, the objects of $P^Q$ over $x$ are the functors $Px\to Qx$,
while the arrows of $P^Q$ over $f$ are the functors in $Pf\to Qf$ over $\2$.

Discrete fibrations and opfibrations are exponentiable. The following particular case 
(which corresponds to Proposition~\ref{26prop1}) is worth noticing:
\begin{prop}  \label{31prop1}
If $A$ is a df and $D$ a dof over $X$ then the exponential $D^A$ in $\CatX$ is a dof 
(and, by duality, $A^D$ is a df).
\end{prop}
\pf
By the above description of exponentials, the objects over $x$ of $D^A$ are the {\em mappings} $Ax\to Dx$.
Furthermore, arrows over $f$ in $D^A$ are morphisms $Af\to Df$ in $\Cat/\2$. 
But $Af$ and $Df$ have discrete fibers and no non-trivial compositions; they look like these:

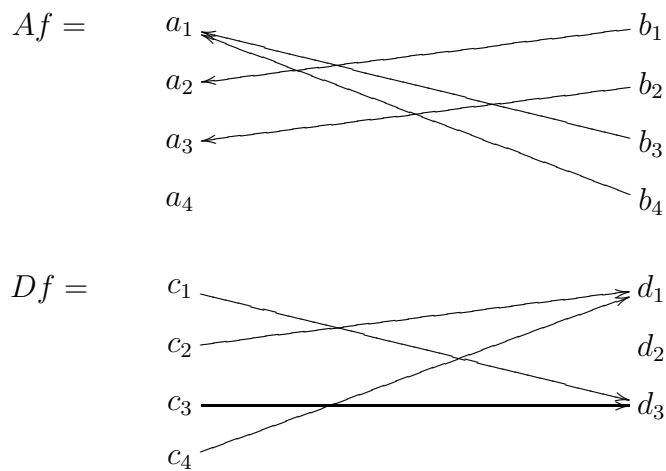
\[ Af = \qq
\xymatrix@R=0.5pc@C=3pc{
{a_1} &&&& \ar[lllld]{b_1} \\
{a_2} &&&& \ar[lllld]{b_2} \\
{a_3} &&&& \ar[lllluu] {b_3} \\
{a_4} &&&& \ar[lllluuu]{b_4} } \\
\] 
\[ Df = \qq
\xymatrix@R=0.5pc@C=3pc{
{c_1}\ar[rrrrdd] &&&& {d_1} \\
{c_2}\ar[rrrru] &&&& {d_2} \\
{c_3}\ar[rrrr] &&&& {d_3} \\
{c_4}\ar[rrrruuu] &&&&  }
\]
and the morphisms $Af\to Df$ over $\2$ are uniquely determined by the mappings $Ax\to Dx$; 
indeed, any such a mapping forces an arrow mapping which, on its own turn, forces $Ay\to Dy$.
For example, if $a_2\mapsto c_3$ then the arrow in $Af$ starting from $b_1$
is forced to be mapped to the arrow in $Df$ starting from $c_3$, so that in particular $b_1\mapsto d_3$, and so on.
Then $D^A$ is indeed a dof. 
\epf
\noindent Explicitly, as a covariant presheaf $D^A$ is the composite
\[ 
\xymatrix@1@C=4pc{ X \ar[r]^-\Delta & X\times X \ar[r]^-{A\times D} & \Set\op\times\Set \ar[r]^-\hom & \Set } 
\]
that is, $(D^A)x = Dx^{Ax}$ and 
\eqa* (D^A)f = Df^{Af}&:&Dx^{Ax}\to Dy^{Ay} \\
                      & & h\mapsto Df\circ h\circ Af 
\eeqa*
The categories over $X$ that are both df's and dof's are called discrete bifibrations (dbf's) and correspond to 
presheaves or to covariant presheaves that act by bijections; 
the direct link between these associates to $A:X\op\to\Set$ the $A':X\to\Set$ defined by $A'x = Ax$ and $A'f = (Af)\inv$. 
In particular, if $X$ is a grupoid any df or dof is a bifibration, and there is an isomorphism $\Set^{X\op} \cong \Set^X$
defined by $A'f = A(f\inv)$.
\begin{corol}  \label{31prop2}
If $B$ and $C$ are discrete bifibrations over $X$ then also the exponential $B^C$ in $\CatX$ is a discrete bifibration.
In particular, if $X$ is a grupoid, exponentials in 
\[ \Set^{X\op}\cong\Set^X\simeq\df X = \dof X \] are computed as in $\CatX$.
For example, if $B,C : X\to\Set$, their exponential in $\Set^X$ is given by $(B^C)x = Bx^{Cx}$ and
\eqa* (B^C)f = Bf^{Cf\inv}&:&Bx^{Cx}\to By^{Cy} \\
                        & & h\mapsto Bf\circ h\circ Cf\inv 
\eeqa*	
\end{corol}
\epf
Among the dbf's there are the constant ones: $\disc S$ with $S\in\Set$.
\begin{corol}  \label{31prop3}
If $A$ is a df and $D$ a dof over $X$ then the exponential $(\disc S)^A$ in $\CatX$ is a dof and $(\disc S)^D$ is a df.
Explicitly, as a covariant presheaf $(\disc S)^A = \Set(A-,S)\,$, that is: $(\disc S)^A x = S^{Ax}$ and 
\eqa* (\disc S)^A f = S^{Af}&:&S^{Ax}\to S^{Ay} \\
                            & & h\mapsto h\circ Af 
\eeqa*                            
\end{corol}
\epf

\subsection{Strong dinaturality}
\label{3sub2}

So as the concept of natural transformation between presheaves $A,B:X\op\to\Set$ is subsumed by that of morphism between 
the corresponding df's over $X$:
\[ \Nat(A,B) \cong  \CatX(A,B) \]
we now show how also the strong dinatural transformations (also known as \v Barr dinatural", see~\cite{par98} and~\cite{din}) 
between set-functors $H,K:X\op\times X\to\Set$ can be seen as morphisms between suitably constructed categories over $X$:
\[ \Dinat^*(H,K) \cong  \CatX(H,K) \] 
giving them a very natural status. 

Given a functor $H:X\op\times X \to \Set$, we define a category over $X$ as follows:
\begin{itemize}
\item
the objects over $x\in X$ are the elements of $H(x,x)$;
\item
given $f:x\to y$ in $X$, there is at most one arrow from $a\in H(x,x)$ to $b\in H(y,y)$ over $f$,
and this is the case iff $H(x,f)a = H(f,y)b \in H(x,y)$. So we can name this arrow also by $f$. 
\end{itemize}
We only have to check composition: if $f:a\to b$ and $g:b\to c$ then $H(x,f)a = H(f,y)b$ and $H(y,g)b = H(g,z)c$. So
\eqa* H(x,gf)a = H(x,g)(H(x,f)a) = H(x,g)(H(f,y)b) = H(f,g)b \\ 
               = H(f,z)(H(y,g)b) = H(f,z)(H(g,z)c) = H(gf,z)c 
\eeqa*
As in the case of presheaves and df's, we often refer to this category as \v $H$ over $X$" (or simply as $H$).
\begin{examples}  \label{32ex1}
The following particular cases are worth noting:
\begin{enumerate}
\item
If $H$ is constant, then $H$ over $X$ is the corresponding constant dbf.
\item
If $H$ is dummy in one variable, that is if it factors through a projection, then we find again
the df or dof of section~\ref{3sub1}.
\item
If $X$ is a grupoid, $H$ over $X$ is the bifibration corresponding to 
\[ \xymatrix@1{X\ar[rr]^-{\Delta_X} && X\times X \ar[rrr]^-{(-)\inv\times X} &&& X\op\times X \ar[rr]^-H && \Set} \]
\item
Given $A:X\op\to\Set$ and $D:X\to\Set$, let $H(x,y)=Ax\times Dy$, that is $H$ is the composite 
\[ \xymatrix@1{X\op\times X \ar[rr]^-{A\times D} && \Set\times\Set \ar[rr]^-\times && \Set} \]
Then $H$ over $X$ is the product of $A$ and $D$ over $X$ (that is, of the df $A$ and the dof $D$ in $\CatX$).
\item
Given $A,B:X\to\Set$, let $H(x,y) = By^{Ax} = \Set(Ax,By):X\op\times X\to\Set$, that is the composite 
\[ \xymatrix@1{X\op\times X \ar[rr]^-{A\times B} && \Set\op\times\Set \ar[rr]^-\hom && \Set} \]
Then it is easy to see that $H$ over $X$ is the exponential $B^A$ in $\CatX$. 
\item
If $H$ is $\hom_X:X\op\times X\to\Set$ then $H$ over $X$ has the endomorphisms in $X$ as objects, and 
given $h:x\to x$ and $k:y\to y$, $f:x\to x'$ is an arrow $h\to k$ in $H$ iff $f\circ h = k\circ f$. 
So $H$ over $X$ is the endomorphism category $X^\N$, where $\N$ is the monoid of natural numbers.
\item
If $F:X\to X$ is an endofunctor on $X$, and $H$ is $\hom_X(F-,-):X\op\times X\to\Set$ then $H$ over $X$
is the category of \v $F$-algebras" with the forgetful functor to $X$.
(The commutative squares below correspond to the arrows in $H$ in the last three examples.)
\eq
\xymatrix@R=3.5pc@C=3.5pc{
Ax    \ar[r]^{Af}\ar[d]_h    &   Ay \ar[d]^k   \\
Bx    \ar[r]^{Bf}            &   By             }
\qq\qq
\xymatrix@R=3.5pc@C=3.5pc{
x    \ar[r]^f\ar[d]_h    &   y \ar[d]^k   \\
x    \ar[r]^f            &   y             }
\qq\qq
\xymatrix@R=3.5pc@C=3.5pc{
Fx    \ar[r]^{Ff}\ar[d]_h    &   Fy \ar[d]^k   \\
x     \ar[r]^f            &   y             }
\eeq
\end{enumerate}
\begin{remark}
While in the first four examples above, $H$ over $X$ is exponentiable in $\CatX$, this is not true in general. 
E.g., let $X$ be the category with three objects and the following non-identity arrows
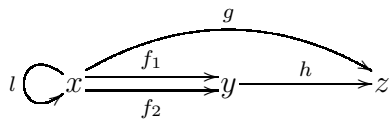
\[ \xymatrix@1{x\ar@(ul,dl)[]_l\ar@<0.5ex>[rr]^{f_1}\ar@<-0.5ex>[rr]_{f_2}\ar@/^1.5pc/@<0.5ex>[rrrr]^g && y \ar[rr]^h && z } \]
with $l^2=x$, $f_1 l = f_2$, $f_2 l = f_1$. 
If $H$ is the endomorphism category of $X$ as in example 6 above, then $g:l\to z$ in $H$ and $g=h f_1$ in $X$,
but this factorization cannot be lifted to $H$, since $f_1$ is not an arrow $l\to y$.  
\end{remark}
\end{examples}
Now, given $H,K:X\op\times X \to \Set$, what is a morphism in $\CatX$ among the corresponding categories over $X$?
First of all, we need to give a family of mappings $\alpha_x : H(x,x) \to K(x,x)$.
If there is a functor over $X$ with such a family as object mapping, it is uniquely determined.
And this is the case if and only if for any $f:x\to y$ in $X$, whenever $f:a\to b$ is an arrow in $H$, 
$f:\alpha_x a \to \alpha_x b$ is an arrow in $K$ too. 
Explicitly, the $\alpha_x$ define a morphism $H\to K$ in $\CatX$ if and only if the following diagram 
commutes for any $f:x\to y$ in $X$: 
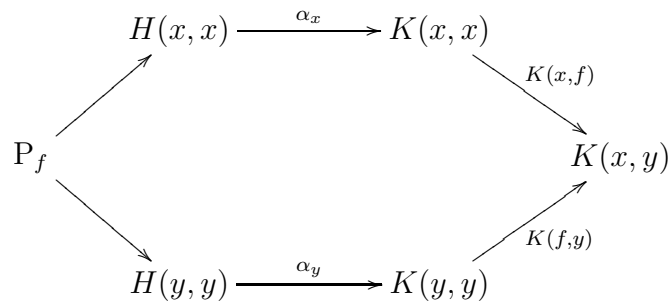
\[
\xymatrix@R=2.5pc{
                  & H(x,x)\ar[rr]^{\alpha_x} && K(x,x)\ar[dr]^{K(x,f)} &  \\ 
\pb_f\ar[ur]\ar[dr] &                          &&                        & K(x,y)  \\ 
                  & H(y,y)\ar[rr]^{\alpha_y} && K(y,y)\ar[ur]_{K(f,y)} &  \\  }
\]
where the pullback $\pb_f$
\eq \label{32eq1}
\xymatrix@R=3.5pc@C=3.5pc{
\pb_f    \ar[r]\ar[d]    & H(x,x) \ar[d]^{H(x,f)} \\
H(y,y) \ar[r]^{H(f,y)} & H(x,y) }
\eeq
expresses the pairs of objects $\langle a,b \rangle$ of $H$ such that $f:a\to b$ is an arrow in $H$.
So morphisms $H\to K$ over $X$ correspond to the strong dinatural transformations $H\to K$. 
Of course, if
\eq \label{32eq2}
\xymatrix@R=3.5pc@C=3.5pc{
H(y,x) \ar[r]^{H(f,x)}\ar[d]_{H(y,f)} & H(x,x) \ar[d]^{H(x,f)} \\
H(y,y) \ar[r]^{H(f,y)} & H(x,y) }
\eeq
is a pullback for any $f:x\to y$, then dinatural transformations with domain $H$ are already strongly dinatural.
Here are some typical cases when this is the case for any codomain category:
\begin{prop}  \label{32prop1}
Let $H:X\op\times X \to \Set$. Each of the following conditions implies that~{\rm (\ref{32eq2})} is a pullback
for any $f:x\to y$, and so also that
\[ \Dinat(H,K) \cong \Dinat^*(H,K) \cong  \CatX(H,K) \]
for any $K:X\op\times X \to \Set\,$: 
\begin{enumerate}
\item
$H$ is constant.
\item
$H$ is dummy in one variable.
\item
$X$ is a grupoid.
\item
$H(x,y)=Ax\times Dy\q$ {\rm (see Examples~\ref{32ex1})}.
\end{enumerate}
\end{prop} 
\pf
For example, let's check the last case (the other ones are simpler):
\eq \label{32eq3}
\xymatrix{
T \ar@/^1.5pc/[drrr]^{\tau_x} \ar@{-->}[dr]^\tau \ar@/_1.5pc/[dddr]_{\tau_y} \\
& Ay\times Dx \ar[rr]^{Af\times Dx}\ar[dd]|{Ay\times Df} && Ax\times Dx \ar[dd]|{Ax\times Df} \\
&                                                        && \\
& Ay\times Dy \ar[rr]^{Af\times Dy}                      && Ax\times Dy }
\eeq
Let's denote by $\ppi^{x,x}{}\!\!_i$ the projections of $Ax\times Dx$, and similarly the other projections.
Given $\tau_x$ and $\tau_y$ which make the diagram~(\ref{32eq3}) commute, if $\tau$ exists we have
\[ 
\ppi^{y,x}\!\!_2\circ\tau = Dx\circ\ppi^{y,x}\!\!_2\circ\tau = \ppi^{x,x}\!\!_2\circ (Af\times Dx)\circ\tau = \ppi^{x,x}\!\!_2\circ\tau_x
\]
and similarly, $\ppi^{y,x}\!\!_1\circ\tau = \ppi^{y,y}\!\!_1\circ\tau_y$
so that $\tau$ must have the form $\tau = \langle \ppi^{y,y}\!\!_1\circ\tau_y , \ppi^{x,x}\!\!_2\circ\tau_x \rangle$.
Now to check e.g. that $(Af\times Dx)\circ\tau = \tau_x$ we compose with projections:
\[ \ppi^{x,x}\!\!_2\circ (Af\times Dx)\circ\langle \ppi^{y,y}\!\!_1\circ\tau_y , \ppi^{x,x}\!\!_2\circ\tau_x \rangle
= Dx\circ\ppi^{y,x}\!\!_2\circ\langle \ppi^{y,y}\!\!_1\circ\tau_y , \ppi^{x,x}\!\!_2\circ\tau_x \rangle = \ppi^{x,x}\!\!_2\circ\tau_x \]
\eqa* 
\ppi^{x,x}\!\!_1\circ (Af\times Dx)\circ\langle \ppi^{y,y}\!\!_1\circ\tau_y , \ppi^{x,x}\!\!_2\circ\tau_x \rangle
= Af\circ\ppi^{y,x}\!\!_1\circ\langle \ppi^{y,y}\!\!_1\circ\tau_y , \ppi^{x,x}\!\!_2\circ\tau_x \rangle = Af\circ\ppi^{y,y}\!\!_1\circ\tau_y \\
= \ppi^{x,y}\!\!_1\circ (Af\times Dy)\circ\tau_y = \ppi^{x,y}\!\!_1\circ (Ax\times Df)\circ\tau_x = \ppi^{x,x}\!\!_1\circ\tau_x                      
\eeqa*
\epf
\noindent The third case of the above proposition, that is when $X$ is a grupoid, can be reversed in the following sense:
\begin{prop}
If the~{\rm (\ref{32eq2})} are pullbacks for any $H$ and $f$ then $X$ is a grupoid. 
\end{prop}
\pf
Let $H = \hom_X:X\op\times X\to\Set$; the pair $\langle \id_x , \id_y \rangle$ is in the pullback~(\ref{32eq1}) for any $f:x\to y$. 
Then, if~(\ref{32eq2}) is such a pullback, there is an arrow $f':y\to x$ in $\hom_X(y,x)$ such that 
\[ f'\circ f = \hom_X(f,x)f' = \id_x \qv f\circ f' = \hom_X(y,f)f' = \id_y \]
\epf
\noindent In section~\ref{3sub5} we shall apply some of the results illustrated to the calculus of ends and coends.

\subsection{Components of graphs and categories}
\label{3sub3}

In order to properly generalize the \v meets" operator of section~\ref{2sub4} to the set-valued context, we need to consider
categories $\A$ that have components, that is there exist functors
\[
\comp \adj \disc \adj \pt : \A\to\Set
\]
called the \v components", the \v discrete" and the \v points" functors.
The objects of $\A$ isomorphic to those of form $\disc S$ (for a set $S\in\Set$) are said to be discrete or constant, 
while elements of $\pt A$ are called points and those of $\comp A$  components of the object $A\in \A$.
If $\A$ is a category with components, then it is such in an essentially unique way, since in this case
$\disc 1$ is terminal in $\A$, $\disc S$ is the copower $S\cdot 1$ (that is the sum of $S$ copies of $1\in \A$) 
and the points functor is represented by $1\in \A$.
 
Among categories with components, of pivotal importance is the category of (irreflexive) graphs:
\[
\comp \adj \disc \adj \pt : \Grph\to\Set
\]
Since the terminal graph is the loop, $\pt X$ is the set of the loops of the graph $X$, while the discrete graphs 
are sums of loops. 

Our attitude in the present paper is of considering the category $\Grph$ as primitive, rather than as a presheaf category:
consider a graph as made up of nodes and arrows between them;
arrows are immaterial, so that they can cross each other, and yet each one bounds his domain and codomain nodes.
The set $\comp X$ of the components of a graph $X$ is what you get by shrinking its arrows to a zero length.
The adjunction $\comp\adj\disc$ expresses the fact that in order to map a graph to a discrete graph, 
one has to perform a similar process, before eventually doing further identifications 
(now the arrows are shrinked to loops, but the effect on the nodes is the same).
Usually one says: \v the components of $X$ are the equivalence classes in $X/\!\!\sim$ of the least equivalence relation $\sim$ 
on the nodes of $X$ which contains the relation $R$ defined by $x R y$ iff there exists an arrow $f:x\to y$ in $X$".
This follows from the fact that, since a discrete graph is in fact an equivalence relation, $\disc:\Set\to\Grph$
factors through $\Set\to\Eq\to\Rel\to\Grph$ and its left adjoint $\comp:\Grph\to\Set$ also factors through
$\Grph\to\Rel$ (\v there exists an arrow"), $\Rel\to\Eq$ (\v the least equivalence relation", which in turn could be
factorized explicitly though symmetrical relations: \v $x\sim y$ iff there is a sequence $x_0=x,x_1,\dots,x_n=y$ 
such that $x_i R x_{i+1}$ or $x_{i+1} R x_i$") 
and $\Eq\to\Set$ (\v the equivalence classes"). Of course, other factorizations are possible, such as that through categories with
involution (reversible categories):$\q\comp:\Grph\to\Cat^{\bf rev}\to\Eq\to\Set$ (see~\cite{hig71}); then one would say
\v the components of $X$ are the equivalence classes of the equivalence relation
on the nodes of $X$ given by $x\sim y$ iff there exists an undirected path $x\to y$ in $X$".

While these factorizations of the left adjoint may be useful to compute it (e.g. to prove that two nodes 
\v are in the same component", that is are mapped to the same element by the unit $\eta:X\to\disc\comp X$ of the adjunction), 
from a conceptual point of view they conceal, rather than explain, the simplicity of the above sketched idea 
of the direct adjunction. (On the contrary, one can see the above categories as categories of special graphs, or graphs with 
structure, so that their own components functors are subsumed by that on graphs.) 

Of course, one can say that graphs are presheaves on the category $\D = \xymatrix@1{\bullet \ar@<0.5ex>[r]\ar@<-0.5ex>[r] & \bullet}$,
so that components are given by the colimit (coequalizer) functor $\Set^{\D\op}\to\Set$. But concretely
colimits in $\Set$ are given by the components of a graph (of the corresponding discrete fibration, see section~\ref{3sub4}).
So we would reduce the components of a graph to those of a more involved graph (the \v subdivision graph" of Kan).  

If $X$ is a graph, $\Grph/X$ also has components: if $\pi:P\to X$ is a graph over $X$, its points are the sections of $\pi$,
while its components are those of the total graph $P$. 

The category $\Cat$ of categories also has components:
\[
\comp \adj \disc \adj \pt : \Cat\to\Set 
\]
the points being the objects of the category, while the components are those of the underlying graph.
Similarly we have
\[ \comp \adj \disc \adj \pt : \CatX\to\Set \]
components being given by those of the total category.
The discrete objects $\disc S$ of $\CatX$ are those discrete bifibrations that correspond to 
the constant presheaves with value $S\in\Set$.

\subsection{Limits and colimits of set functors}
\label{3sub4}

Of course, any presheaf category $\Set^{X\op}$ has components
\[
\comp \adj \disc \adj \pt : \Set^{X\op}\to\Set
\]
Usually $\pt$ and $\comp$ are known as the limit and the colimit functors, while as $\disc$ one can use the diagonal functor:
\[
\Colim \adj \Delta \adj \Lim : \Set^{X\op}\to\Set
\]
The important fact is that, since $\disc_{\CatX}$ factors (up to an isomorphism) through $\disc_{\Set^{X\op}}$, 
\eq \label{34eq0}
\xymatrix@R=4pc@C=4pc{
{\Set^{X\op}} \ar@{_{(}->}[d] & {\Set}\ar[l]^-\Delta\ar@2{-}[d] \\ 
{\CatX}                       & {\Set}\ar[l]_-\disc }
\eeq
these functors can be obtained by restricting the corresponding ones of $\CatX$:
\eq \label{34eq1}
\xymatrix{
{\Set^{X\op}} \ar[rr]^\Lim\ar@{_{(}->}[dd] && {\Set}\ar@2{-}[dd] \\ \\
{\CatX}       \ar[rr]^\pt                        && {\Set} }
\qq\qq
\xymatrix{
{\Set^{X\op}} \ar[rr]^\Colim\ar@{_{(}->}[dd] && {\Set}\ar@2{-}[dd] \\ \\
{\CatX}       \ar[rr]^\comp                        && {\Set} }
\eeq
so that the limit of $A:X\to\Set$ is the set of sections of $A$ over $X$, while its colimit is the set components of $A$ over $X$ 
(that is, the components of the graph of the category of elements of $A$).
More explicitly, we have the following chain of natural isomorphisms:
\[
\begin{array}{c}
\Cone(S,A) \\ \hline
\Nat(\Delta S,A) \\ \hline
\CatX(\disc S,A) \\ \hline
\Set(S,\pt A) 
\end{array}
\qq\qq
\begin{array}{c}
\Cone(A,S) \\ \hline
\Nat(A,\Delta S) \\ \hline
\CatX(A,\disc S) \\ \hline
\Set(\comp A,S) 
\end{array}
\]
thus $\pt A = \Lim A$ and $\comp A = \Colim A$.
\begin{remark}  \label{34rmk1}
We have just seen that colimits of set functors $A:X\op\to\Set$ are given by the components of the df $A$ over $X$.
Conversely, the components of any category $P\in\Cat$ can be expressed as the colimit
of a functor $X\op\to\Set$, for {\em any} df $\, P\to X$ with total category $P$.
In particular $P$ is a df $\, P\tto^\id P$ on itself, corresponding to the functor $\Delta 1 : P\op\to\Set$, 
the terminal object in $\Set^{P\op}$, so that 
\[ \comp P \,\cong\, \Colim \Delta 1 \] 
\end{remark}

\subsection{Ends and coends}
\label{3sub5}

As remarked in Proposition~\ref{32prop1}, if the domain is constant strong dinaturality reduces to standard dinaturality:
\[ \Dinat(\Delta S,H) \cong \Dinat^*(\Delta S,H) \]
So we have
\[
\begin{array}{c}
\Wedge(S,H) \\ \hline
\Dinat(\Delta S,H) \\ \hline
\Dinat^*(\Delta S,H) \\ \hline
\CatX(\disc S,H) \\ \hline
\Set(S,\pt H) 
\end{array}
\]
and we conclude
\[ \pt H = \End H \]
For coends in general the situation is slightly different:
\[
\begin{array}{c}
\Wedge^*(H,S) \\ \hline
\Dinat^*(H,\Delta S) \\ \hline
\CatX(H,\disc S) \\ \hline
\Set(\comp H,S) 
\end{array}
\]
so that we can only conclude that
\[ \comp H = \Coend^* H \]
where $\Coend^* H$ is a representing object for $\Wedge^*(H,S)$.
But, as shown in the example below, in general $\Coend^* H \ne \Coend H$.
\begin{examples} 
Let us test the above concepts in the case $H = \hom_X$.
Recall that $\hom_X$ over $X$ is the endomorphism category of $X$ (see the sixth of Examples~\ref{32ex1}).
\begin{itemize}
\item
$\End\,\hom_X = \End^*\hom_X = \pt\hom_X$ is the set of its sections: families of endomorphisms $\, e_x:x\to x,\, x\in X$, 
such that $f\circ e_x = e_y\circ f$ for any $f:x\to y$ in $X$. This is the well-known \v center" of $X$. 
\item
$\Coend^*\hom_X = \comp\hom_X$, wherein two endomorphisms $e_x$ and $e_y$ are identified whenever $f\circ e_x = e_y\circ f$, 
with $f:x\to y$ in $X$.
\item
$\Coend\,\hom_X$ is obtained by identifying the two compositions $f\circ g$ and $g\circ f$ of two cycles of arrows in $X$: 
\[ \xymatrix{ x\ar@/^/[rr]^f && y\ar@/^/[ll]^g } \]
\end{itemize}
That $\Coend\,\hom_X \ne \Coend^*\hom_X$ is then seen already in the simple case $X=\2$
(denote by $f:x\to y$ its non-identity arrow):
$\id_x$ and $\id_y$ are identified in the strong coend by $f$, while they are not identified in the standard coend
since there are no non-trivial two-cycles. So $\Coend^*\hom_X$ is a terminal set, while $\Coend\,\hom_X$ has two elements.
A similar conclusion holds for the commutative monoid $\{ 1,e \}$ with $e$ idempotent, while in general
if $X$ is a commutative monoid $\Coend\,\hom_X \cong \End\,\hom_X$ is the set of arrows of $X$. 
If $X$ is a grupoid, we know from Proposition~\ref{32prop1} that $\Coend\,\hom_X = \Coend^*\hom_X$; 
in fact they give the conjugacy classes of $X$.
\end{examples}
\begin{remark}  \label{35rmk1}
Following the fifth of Examples~\ref{32ex1}, we have
\[ \Nat(A,B) = \CatX(A,B) = \pt(B^A) = \End\,\Set(A-,B-) \]
that is we find again the well-known expression of the set of natural transformations as an end. 
\end{remark}
The above considerations are summarized by the following commutative diagrams
(which subsume diagrams~(\ref{34eq1}) about limits and colimits),
where $\Set^{X\op\times X}_*$ is the category whose arrows are strong dinatural transformations:
\eq 
\xymatrix{
{\Set^{X\op\times X}_*} \ar[rr]^-{\End}\ar@{_{(}->}@<1ex>[dd] && {\Set}\ar@2{-}[dd] \\ \\
{\CatX}               \ar[rr]^-{\pt}                        && {\Set} }
\qq\qq
\xymatrix{
{\Set^{X\op\times X}_*} \ar[rr]^-{\Coend^*}\ar@{_{(}->}@<1ex>[dd] && {\Set}\ar@2{-}[dd] \\ \\
{\CatX}               \ar[rr]^-{\comp}                          && {\Set} }
\eeq

\section{The tensor functor and complements of categories over a base}
\label{complements}

Following Section~\ref{posets}, we now generalize some aspects of the two-valued duality to the set-valued context:
the \v meets" predicate becomes a \v tensor" functor, $\ten:\A\times\A\to\Set$, which allows to define the \v complement" of $A\in\A$
as the right adjoint to $\ten(A,-)$.
Once more, it appears that $\CatX$ encloses and clarifies well-known concepts: in this case, the tensor functor on $\CatX$ 
extends the classical tensor product $\otimes:\Set^{X\op}\times\Set^X\to\Set$, while each complement functor $\neg A:\Set\to\CatX$
\v extends" the right adjoints of $A\otimes -:\Set^X\to\Set$.

\subsection{The tensor functor and complements}
\label{4sub1}

If $\A$ has components and products we can naturally generalize the \v meets" predicate of Section~\ref{posets} obtaining the functor: 
\[ \ten_\A := \comp (-\times -) : \A\times \A\to\Set \] 
We may call functors isomorphic to some $\ten(A,-):\A\to\Set$ \v t-representable" functors.
It is worth stressing that non isomorphic objects may t-represent isomorphic functors; 
an example is given by two graphs in $\Grph$ with the same thin reflection or 
by any two {\em reflexive} graphs with the same number of components. 
On the other hand, in the category of presheaves on a grupoid it is true that
\[ \ten(A,-) \,\cong\, \ten(B,-) \iimp A\cong B \]
If $\A$ is also (cartesian) closed, then the dual or complementary roles played by the functors $\ten_\A$ and $\hom_\A$ 
are more evident, as summarized in the following proposition: 
\begin{prop} \label{41prop1}
Let $\A$ be a cartesian closed category with components, then
\begin{eqnarray}
\ten(A,B) := \comp (A\times B) \nn \\
\hom(A,B) \cong \pt (B^A) \nn \\ 
\comp \cong \ten(1,-) \nn \\
\pt \cong \hom(1,-) \nn \\
\label{41eq1} 
\ten(A,-) \adj (\disc -)^A \\ 
\label{41eq2}
(\disc -)\times A \adj \hom(A,-) 
\end{eqnarray}
the last two being adjunctions with parameter $A\in\A$.
\end{prop}
\epf
In general, the copower functor $-\cdot A:\Set\to \A$, left adjoint to the representable 
$\hom(A,-)$, may exist without being given by the formula displayed in equation~(\ref{41eq2}) (as in the case of categories
with a zero object, on which the functors $\comp$, $\disc$ and $\pt$ become constant). 
On the other hand, we will shortly show that if $\ten(A,-)$ has a right adjoint then it has the form displayed in~(\ref{41eq1}).
\begin{definition} \label{41def1}
If there is a universal arrow $\ten(A,B)\to S$ from $\ten(A,-)$ to the set $S$, we denote $B$ by $(\neg A)S$
and call it \v the $S$-complement of $A$". 
In particular, if $(\neg A)S$ exists for any set $S$, we say that $A$ \v has complement" and call functor $\neg A:\Set\to\A$, 
right adjoint to $\ten(A,-)$, \v the complement of $A$". 
\end{definition}
 
\begin{prop} \label{41prop2}
If $\A$ has products and components, an object $A\in\A$ has an $S$-complement iff the exponential $(\disc S)^A$ exists
(that is, iff there is a universal arrow from $A\times -$ to the set $\disc S$).
If this is the case, 
\[ (\neg A)S \cong (\disc S)^A \] 
In particular, $A$ has complement iff the exponential $(\disc S)^A$ exists for any set $S$. 
If this is the case (that is, if $A$ is \v exponentiable for discrete basis") then 
\[ \neg A \cong (\disc -)^A \] 
\end{prop}
\pf
If $(\disc S)^A$ exists, its universality with respect to $\ten(A,-)$ and $S$
is given by the following chain of bijections, natural in $B\in\A$: 
\[
\begin{array}{c}
\ten(A,B) \to S \\ \hline 
\comp(A\times B) \to S \\ \hline 
A\times B \to \disc S \\ \hline
B\to (\disc S)^A 
\end{array}
\]
In the other direction, if $(\neg A)S$ exists, we have
\[
\begin{array}{c}
B\to (\neg A)S \\ \hline
\ten(A,B) \to S \\ \hline 
\comp(A\times B) \to S \\ \hline 
A\times B \to \disc S 
\end{array}
\]
which gives the universality of it with respect to $A\times -$ and $\disc S$.
\epf
\noindent In general, the adjunction with parameter $A\in\A$:
\[ \ten(A,-) \adj \neg A = (\disc -)^A \]
is only partially defined.

\subsection{Atoms}
\label{4sub2}

Let $\A$ be a category with products and components.
Following section~\ref{2sub5}, we say that an object $A\in\A$ is a \v strong atom" if there are bijections
\[ \ten(A,B) \,\cong\, \hom(A,B) \]
natural in $B\in\A$.
For example, if $\A=\Grph$ is the category of graphs, the dot graph (one object and no arrows) is a strong atom.
The set-valued correspective of Proposition~\ref{25prop1} (to which we refer for the proof) is the following:
\begin{prop}  \label{42prop1}
If a strong atom $A$ of $\A$ has complement then it is a absolutely presentable object (that is, $\hom(A,-)$
preserves all colimits).
In particular, if $\A$ is cartesian closed then its strong atoms are absolutely presentable objects.
\end{prop}
\epf
\noindent Often the above condition is too strong to be interesting.
But it can be relativized by restricting the class of objects against which to test atomicity.

In our main concern, $\A=\CatX$, we say that a category $P$ over $X$ is a \v left atom" 
if there is a natural bijection as above
\[ \ten(P,A) \,\cong\, \hom(P,A) \]
but restricted to the $A\in\df X$. Dually, $P$ is a \v right atom" if 
\[ \ten(P,D) \,\cong\, \hom(P,D) \]
natural in $D\in\dof X$.
We say that $P\in\CatX$ is an atom if it is both a left and a right atom.
We shall study thoroughly atoms in section~\ref{6sub2}, where we shall also see that
the distinction between left and right atoms is illusory.
For now we just record
\begin{prop}  \label{42prop2}
Any object $x$ of the base category $X$ is an atom: for any df or dof $A$
\[ \ten(x,A) \cong \hom(x,A) \cong Ax \]
\end{prop}
\pf
In general $\ten(x,P)=\comp Px$ and $\hom(x,P)=\pt Px$ are the components and the objects of the fiber category
over $x$. But a category $\A$ is discrete iff $\comp\A \cong \pt\A$, and so if $P$ has discrete fibers 
then $\comp Px \cong \pt Px$.
\epf

\subsection{The tensor functor on $\CatX$ and tensor products}
\label{4sub3}

We now show that $\ten_{\CatX}:\CatX\times\CatX\to\Set$ extends the usual tensor product of set functors.
Recall that given $A:X\op\to\Set$ and $D:X\to\Set$, if $H(x,y)=Ax\times Dy$ like in Examples~\ref{32ex1}, then
\eq  \label{43eq1}
A\otimes D := \Coend H   
\eeq
is called the tensor product of $A$ and $D$.
We have seen that $H$ over $X$ is the product of $A$ and $D$ over $X$, and that in this case 
\[ \Wedge(H,S) = \Wedge^*(H,S) \]
(see Proposition~\ref{32prop1}). So we have 
\[ \Coend H = \Coend^* H = \comp H = \comp(A\times B) = \ten_{\CatX}(A,B) \]
that is, the following diagram commutes:
\eq   \label{43eq2}
\xymatrix{
{\Set^{X\op}\!\!\times\Set^X} \ar[rr]^-\otimes\ar@{_{(}->}[dd] && {\Set}\ar@2{-}[dd] \\ \\
{\CatX\times\CatX}            \ar[rr]^-\ten                          && {\Set} }
\eeq
So the tensor product is reduced directly to the components of a category over $X$.
Usually it is reduced to the components of the same (total) category but seen as a dof over a df over $X$:
a colimit of the set functor 
\[ A \ttto^{\pi_A} X \ttto^D \Set \]
on the category of elements of $A$. The two constructions are clearly equivalent, but the former is more direct and symmetrical. 
The other possibility would be the reduction of the coend~(\ref{43eq1}) to a colimit via the \v subdivision category" or the
\v twisted arrow category" (see~\cite{mac71}). But so doing would mean to consider the components of a more involved 
total category (whose objects are given by all the $Ax\times Dy$, not only those with $x = y$).

Observe that \v dually" the diagrams
\[ \xymatrix{
{\Set^{X\op}\!\!\times\Set^{X\op}} \ar[rr]^-\Nat\ar@{_{(}->}[dd] && {\Set}\ar@2{-}[dd] \\ \\
{\CatX\times\CatX}                 \ar[rr]^-\hom                       && {\Set} }
\qq\qq
\xymatrix{
{\Set^X\!\!\times\Set^X} \ar[rr]^-\Nat\ar@{_{(}->}[dd] && {\Set}\ar@2{-}[dd] \\ \\
{\CatX\times\CatX}       \ar[rr]^-\hom                       && {\Set} }
\]
express the fullness of the inclusions~(\ref{31eq2}).

\subsection{Negation and complements in $\CatX$}
\label{4sub4}

Now we turn to complements in $\CatX$ and will see in what sense the dualities of Section~\ref{posets}
can be generalized to the set-valued context.
We have seen in section~\ref{2sub6} that the truth value category $\2$ is very peculiar in that
the functor $\q\no := -\imp\false:\2\to\2\op\q$ is an isomorphism, which we usually call \v negation"
and which in turn induces a duality between lower and upper parts of a poset $\neg:\df X\to(\dof X)\op$.
The point to be made presently is that for the category of sets there is not a single dualizing object
playing the role of \v false"; rather, we must consider the whole bunch of functors
\[ \no_S := S^- : \Set\to\Set\op \]
as the \v negation" in $\Set$, which in turn induce a bunch of negation operators
\[ \neg_S:\df X\to(\dof X)\op \]
These are not isomorphisms, but negation is classical in that a presheaf $A$ is recoverable from the
bunch (in fact a functor on $\Set$) of covariant presheaves $\neg_S A$.
Furthermore, once again, the negation operator extends, at least partially, to the environment $\CatX$
making it possible to derive the formulas for the reflection in discrete fibrations: the $S$-negation
of a presheaf is the \v same" as the $S$-complement in $\CatX$ of the corresponding df, as in the two-valued case. 
(In~\cite{law96}, it is proposed a definition of negation in a locally cartesian closed category
or in any hyperdoctrine, that is similar to our definition of complement but is parametrized by
objects over 1, rather than by sets).

Let us begin with complements in the category of sets. 
Since $\Set$ is cartesian closed and
\[ \comp \adj \disc \adj \pt : \Set\to\Set \]
are the identity functors, any set $S$ has the complement functor
\[ (-)^S : \Set\to\Set \]
which we call also the \v negation" of $S$, since we are thinking to $\Set$ as a \v truth-values" category. 
So the $T$-negation of $S$ is simply $\no_T S = T^S$, where $\no_T$ is the functor
\[ \no_T = T^- : \Set\to\Set\op \]
As for the complement operator on $\CatX$,
\[ \neg P = (\disc -)^P \]
in general it is only partially defined. 
Anyway Corollary~\ref{31prop2} can now be rephrased as follows:
\begin{prop}  \label{44prop1}
Any df $A$ on $X$ has a complement 
\[ \neg A = (\disc-)^A:\Set\to\CatX \]
which is valued in dof's (and conversely).
So the $S$-complement of a df on $X$ is a dof (and conversely).  
\end{prop}
\epf
On the other hand, as sketched before, the dof $(\neg A)S$ can also be seen as induced directly by negation in $\Set$:
\[ (\no_S)^{X\op} : \Set^{X\op}\to(\Set\op)^{X\op} \cong (\Set^X)\op \]
gives a functor
\[ \neg_S : \Set^{X\op}\to(\Set^X)\op \] 
Expicitly, for any $A:X\op\to\Set$, $\neg_S A$ is the opposite of $(S^-)\circ A$, that is 
\[ \neg_S A = \Set(A-,S) \]
So $\neg_S A$ corresponds to the exponential $(\disc S)^A$ in $\CatX$ for the df $A$, as described in Corollary~\ref{31prop2},
and we can construct exponentials in $\CatX$ in such a way that the following diagram commutes for all $S\in\Set$:
\eq  \label{44eq1}
\xymatrix{
{\Set^{X\op}} \ar[rr]^-{\neg_S}\ar@{_{(}->}[dd] && {(\Set^X)\op}\ar@{_{(}->}[dd] \\ \\
{\CatX}       \ar[rr]^-{(\disc S)^-}            && {(\CatX)\op} }
\eeq
where the lower arrow is the partially defined functor which takes a category over $X$ in its $S$-complement.
Incorporating the parameter $S$ explicitly in the diagram, we get
\eq
\xymatrix{
{\Set\times\Set^{X\op}} \ar[rr]^-{\neg_-}\ar@{_{(}->}@<-1ex>[dd] && {(\Set^X)\op}\ar@{_{(}->}[dd] \\ \\
{\Set\times\CatX}       \ar[rr]^-{(\disc -)^-}                   && {(\CatX)\op} }
\eeq
so, recalling diagram~(\ref{43eq2}), we conclude that 
\[ A\otimes - \adj \neg_- A : \Set\to\Set^X \]
is an adjunction with parameter $A\in\Set^{X\op}$. 
Indeed, this is the usual closed structure on a fragment of the closed bicategory of distributors.
So, the point to be stressed is that the above \v classical" adjunction extends to the adjunction
\[ \ten(A,-) \adj \neg A : \Set\to\CatX \]
of section~\ref{4sub1}, with parameter $A$ in $\df X$.

\subsection{Properties of the tensor functor and of complements in $\CatX$}
\label{4sub5}

Like in section~\ref{2sub8} we now prove an \v adjunction-like" property of the tensor functor on $\CatX$ 
which allows to derive directly the reflection formula.
\begin{prop} \label{45prop1}
Let $X$ be a category, $A$ a df and $D$ a dof on $X$. Then for any $P\in\CatX$
with reflections $\down P$ in df's and $\up P$ in dof's, there are natural bijections
\[ \ten(P,A) \cong \ten(\up P,A) \qv \ten(P,D) \cong \ten(\down P,D) \]
\end{prop}
\pf
We prove the first one, the other being symmetrical.
For any set $S$
\[
\begin{array}{c}
\ten(P,A) \to S \\ \hline 
P \to (\neg A)S \\ \hline 
\up P \to (\neg A)S \\ \hline 
\ten(\up P,A) \to S
\end{array}
\]
where we have used Proposition~\ref{44prop1}. The result follows by Yoneda.
\epf

The following proposition can be interpreted as the fact that complementation in $\CatX$ is classical, 
when restricted to df's and dof's (compare it with Proposition~\ref{28prop2}); 
in particular, it is possible to \v recover" a df (or presheaf) from its negation:
\begin{prop}  \label{45prop2}
If $A$ and $B$ are both df's (or both dof's), then there is a natural bijection
\eq   \label{45eq1}
\hom(A,B) \,\cong\, \hom(\neg B,\neg A)
\eeq
Briefly, the \v strong contraposition law" holds.
Furthermore, df's are properly analyzed by dof's (and conversely) via the tensor functor: 
there is a natural bijection
\[ \hom(A,B) \cong \Nat(\ten(A,-),\ten(B,-)) \]
\end{prop}
\pf
By section~\ref{4sub4}, the isomorphism~(\ref{45eq1}) may be interpreted equivalently using preshaves instead of df's:
\[ \Nat(A,B) \cong \Nat(\neg B,\neg A) \]
where $\neg B,\neg A:\Set\to\Set^X$ are given by $(\neg A)S = \neg_S A = \Set(A-,S)$. 
By the isomorphism $(\Set^X)^\Set \cong \Set^{X\times\Set}$, $\neg A$ corresponds to 
the functor $\neg A : X\times\Set \to \Set$ given by
\[ \neg A (x,S) = (\neg A S)x = \Set(Ax,S) \]
and natural transformations $\alpha:\neg B\to\neg A:\Set\to\Set^X$ correspond to 
natural transformations $\alpha:\neg B\to\neg A:X\times\Set \to \Set$, that is families of mappings
\[ \alpha_S^x:\Set(Bx,S)\to\Set(Ax,S) \q x\in X\, ,\, S\in\Set \]
such that both the diagrams
\eq  \label{45eq2}
\xymatrix{
{\Set(Bx,S)} \ar[rr]^-{\alpha_S^x}\ar[dd]|{\Set(Bx,f)} && {\Set(Ax,S)}\ar[dd]|{\Set(Ax,f)} \\ \\
{\Set(Bx,T)} \ar[rr]^-{\alpha_T^x}                     && {\Set(Ax,T)}      }
\qq\qq
\xymatrix{
{\Set(Bx,S)} \ar[rr]^-{\alpha_S^x}\ar[dd]|{\Set(Bu,S)} && {\Set(Ax,S)}\ar[dd]|{\Set(Au,S)} \\ \\
{\Set(By,S)} \ar[rr]^-{\alpha_S^y}                     && {\Set(Ay,S)}      }
\eeq
commute for any $u:x\to y$ in $X$ and any $f:S\to T$ in $\Set$. 
By Yoneda, the former of~(\ref{45eq2}) implies that $\alpha_S^x = \Set(h_x,S)$ for a 
unique map $h_x:Ax\to Bx$, so that the latter becomes
\eq   \label{45eq3}
\xymatrix{
{\Set(Bx,S)} \ar[rr]^-{\Set(h_x,S)}\ar[dd]|{\Set(Bu,S)} && {\Set(Ax,S)}\ar[dd]|{\Set(Au,S)} \\ \\
{\Set(By,S)} \ar[rr]^-{\Set(h_y,S)}                     && {\Set(Ay,S)}      }
\eeq  
But if $S$ has more than one element, then for any $h,k:U\to V$ in $\Set$, 
\[ \Set(h,S) = \Set(k,S) : \Set(V,S) \to \Set(U,S) \iimp h = k \]
and so the commutativity of~(\ref{45eq3}) for any $S$ and $u:x\to y$ implies that also
\[
\xymatrix{
Ay \ar[rr]^-{h_y}\ar[dd]_{Au} && By\ar[dd]^{Bu} \\ \\
Ax \ar[rr]^-{h_x}             && Bx             }
\]
commutes for any $u:x\to y$ in $X$, giving a natural transformation $a:A\to B$.
In the other direction, any such $a$ induces $\alpha:\neg B \to \neg A$, and the two processes are 
clearly each other inverses.
The last part of the proposition follows from the fact that
\[ \Nat(\ten(A,-),\ten(B,-)) \cong \Nat(\neg B,\neg A) \]
by a general property of adjunctions with parameters (see e.g.~\cite{bip}).
\epf

\section{The reflection and the coreflection of categories over a base in discrete fibrations}
\label{proof}

We begin by showing that the reflections of the objects of the base category in df's have a very familiar form.
Then we derive the reflection and coreflection formulas for general categories over $X$.
Finally, we verify that the formulas give indeed the desired reflection and coreflection.

\subsection{The reflections of objects and the Yoneda Lemma}
\label{5sub1}

First note that the representable functors on $X$ correspond to slices of $X$:
$X(-,x)$ over $X$ is $X/x$, and dually $X(x,-)$ over $X$ is $x/X$.
Thus for any df $A$, $\,\Nat(X(-,x),A) \cong \CatX(X/x,A)$, and the Yoneda Lemma may then be rephrased as
\[ \hom(X/x,A) \,\cong\, Ax \] 
On the other hand, since also (see section~\ref{4sub2}) \(\, \hom(x,A) \cong Ax\, \), we have
\[ \hom(X/x,A) \,\cong\, \hom(x,A) \]
that is 
\begin{prop}  \label{51prop1}
The reflections of an object-atom $x$ of the base category $X$ in df's and dof's are given by:
\[ \down x = X/x \qv \up x = x/X \]
\end{prop}
\epf
\noindent Since also \( \, \ten(x,A) \cong Ax \,\) (see Proposition~\ref{42prop2}),
by Proposition~\ref{45prop1} we get 
\[ \ten(\up x,A) \,\cong\, Ax \]
that is, recalling diagram~(\ref{43eq2}), we have proved the \v co-Yoneda Lemma": 
\begin{corol}  \label{51prop2}
For any presheaf $A:X\op\to\Set$ and any $x\in X$
\[ X(x,-) \otimes A \,\cong\, Ax \] 
\end{corol}
\epf
\noindent Moreover, if $\Y$ is the functor 
\begin{eqnarray}  \label{51eq1}
\Y:X\to\CatX \\
x\mapsto X/x  \nn
\end{eqnarray}
then the \v addendum" to the Yoneda lemma (see~\cite{mac71}) may be rephrased by saying that the 
equivalence~(\ref{31eq1}) can be obtained as
\begin{eqnarray}   \label{51eq2}
\hom(\Y-,-):\df X \to \Set^{X\op} \\
A\mapsto\hom(\Y-,A):X\op\to\Set \nn
\end{eqnarray}
\v Dually" if we define $\,\Y':X\op\to\CatX\,$ by $\, x\mapsto x/X\,$,
the co-Yoneda lemma may be rephrased by saying that
\begin{eqnarray}   \label{51eq4}
\ten(\Y'-,-):\df X \to \Set^{X\op} \\
A\mapsto\ten(\Y'-,A):X\op\to\Set     \nn 
\end{eqnarray}
is an equivalence. In fact, the functors~(\ref{51eq2}) and~(\ref{51eq4}) are isomorphic (see~\cite{bip}).

\subsection{The coreflection and reflection formulas}
\label{5sub2}

We are now in a position to derive the (coreflection and the) reflection formula in a straightforward manner, 
exactly like in section~\ref{2sub9}, using the object-atoms $x\in X$ (see section~\ref{4sub2}) as test shapes 
both for figures (the elements of $\hom(x,P)$) and for \v t-figures" (the elements of $\ten(x,P)$).

Though in the following propositions we should use the hypothetical form \v if the reflection (coreflection) 
exists then it has the form \dots ", we do not since it will be proved in the next section
that these formulas actually give the desired coreflection and reflection for any category over $X$.

\begin{prop}  \label{52prop1}
The fibers of the coreflections $P\rdown$ and $P\rup$ of a category $P$ over $X$ in df's and dof's have the form
\[ (P\rdown)x = \hom(\down x,P) \qv (P\rup)x = \hom(\up x,P) \]
\end{prop}
\pf
\[
\begin{array}{c}
(P\rdown)x \\ \hline
\hom(x,P\rdown) \\ \hline 
\hom(\down x,P\rdown) \\ \hline 
\hom(\down x,P) 
\end{array}
\qq\qq
\begin{array}{c}
(P\rup)x \\ \hline
\hom(x,P\rup) \\ \hline 
\hom(\up x,P\rup) \\ \hline 
\hom(\up x,P) 
\end{array}
\]
where only the adjunction laws have been used. 
\epf
\begin{prop}  \label{52prop2}
The fibers of the reflections $\down P$ and $\up P$ of a category $P$ over $X$ in df's and dof's have the form
\[ (\down P)x = \ten(\up x,P) \qv (\up P)x = \ten(\down x,P) \]
\end{prop}
\pf 
\[
\begin{array}{c}
(\down P)x \\ \hline
\ten(x,\down P) \\ \hline 
\ten(\up x,\down P) \\ \hline 
\ten(\up x,P) 
\end{array}
\qq\qq
\begin{array}{c}
(\up P)x \\ \hline
\ten(x,\up P) \\ \hline 
\ten(\down x,\up P) \\ \hline 
\ten(\down x,P) 
\end{array}
\]
where we used twice the adjunction-like properties of Proposition~\ref{45prop1}. 
\epf
\noindent It is worth stressing the strong similarity between the derivations of the reflection and coreflection formulas.

To obtain the arrows of, say, $\up P$ (that is, the action of $X$ on its fibers) we should look at the reflections in $\df X$ 
of the arrows in $X$. If $f:x\to y$ is an arrow in $X$, its reflection $\down f$ as a category $f:\2\to X$ over $X$ is $X/x$ again; 
and the domain and codomain inclusions over $X$, $x\to f$ and $y\to f$ become
\[ \id : X/x \to X/x \qv -\circ f : X/x \to X/y \]
So we get the following proposition, of which the Yoneda and co-Yoneda Lemmas are particular cases (see section~\ref{5sub1}):
\begin{corol}  \label{52prop3}
The coreflection $P\rdown$ and the reflection $\down P$ are given as presheaves by
\[ P\rdown\, = \hom(\Y-,P) \qv \down P = \ten(\Y'-,P) \]
\end{corol}
\epf

\subsection{The formulas work}
\label{5sub3}

As it often happens, after having derived the form of a supposed adjoint functor, one has to check that what he
got really does the job. We concentrate on dof's, though by duality one clearly has specular results for df's.

We begin by analyzing morphisms in $\CatX$ to and from dof's.
If the codomain is a dof, such a morphism is easily seen to be determined by its object mapping:
\begin{prop}  \label{53prop1}
Given a category $P$ over $X$ and a dof $D$ on $X$, a mapping $\phi$ over $X$ between the objects of the 
total categories of $P$ and $D$ extends (uniquely) to a functor over $X$ iff  
\[ f(\phi a) = \phi b \]
for any $u:a\to b$ over $f:x\to y$.
\end{prop}
\epf
On the other hand, a morphism in $\CatX$ with a dof as domain has the following form, which is simply
the general form adapted to the arrows of a dof: 
\begin{prop} \label{53prop2}
A functor $\phi$ over $X$ from the dof $D$ to $P$ is given by its object mapping over $X$: 
\[ a \,\mapsto\, \phi a \]
and by its arrow mapping over $X$:
\[ \la a,x\tto^f y \ra\, \mapsto\, \phi(a,x\tto^f y) \] 
(with $a\in Dx$) which is a graph morphism and respects identities and compositions:
\begin{enumerate}
\item
$\phi(a,x\tto^f y):\phi a\to \phi(fa)$
\item
$\phi(a,x\tto^\id x) = \id_{\phi a}\,\,\,$ and $\,\,\,\,\phi(a,x\tto^{gf}z) = \phi(fa,y\tto^g z)\circ \phi(a,x\tto^f y)$
\end{enumerate}
\end{prop}
\epf
Now we describe more concretely $\up P$ and $P\rup$. 
Observing that $P/x = X/x\times P$ and so $\ten(\down x,P) = \comp(P/x)$,
we get the following more familiar version of Proposition~\ref{52prop2} and Corollary~\ref{52prop3}
(see~\cite{law73} and~\cite{par73}):
\begin{corol}  \label{53prop3}
The reflection of a category $P$ over $X$ in discrete opfibrations is given by
\[ (\up P)x = \comp(P/x) \]
for any $x\in X$, and  
\[ f[\la a , \pi a\tto^h x\ra] = [\la a , \pi a\tto^{fh}y\ra] \in \comp(P/y) \]
for any $f:x\to y$ in $X$ and $[\la a , \pi a\tto^h x\ra] \in \comp(P/x)$. 
\end{corol}
\epf
\noindent As for the coreflection, from Proposition~\ref{52prop1}, Corollary~\ref{52prop3} 
and Proposition~\ref{53prop2} we straightforwardly get:
\begin{prop} \label{53prop4}
The coreflection of a category $P$ over $X$ in discrete opfibrations is given by $(P\rup)x = \CatX(\up x,P)$
for any $x\in X$; hence each $\xi\in(P\rup)x$ is a mapping 
\[ x\tto^h y \,\,\mapsto\,\xi(x\tto^h y)\, \in Py \] 
for any arrow in $X$ with domain $x$, and a mapping 
\[ \la x\tto^h y,y\tto^k z\ra\, \mapsto\, \xi(x\tto^h y,y\tto^k z) \] 
for any two consecutive arrows in $X$ with domain $x$, such that $\xi(x\tto^h y,y\tto^k z)$ is an arrow in $P$ over $k$ and
\begin{enumerate}
\item
$\xi(x\tto^h y,y\tto^k z):\xi(x\tto^h y)\to\xi(x\tto^{kh} z)$
\item
\(
\xi(x\tto^h y,y\tto^\id y) = \id_{\xi h} \\
\xi(x\tto^h y,y\tto^{lk}w) = \xi(x\tto^{kh}z,z\tto^l w)\circ\xi(x\tto^h y,y\tto^k z)
\)
\end{enumerate}
Furthermore, for any $f:x\to y$ in $X$, $f\xi :\,\up y \to P$ is given by
\[ (f\xi)(y\tto^h z) = \xi(x\tto^{hf}z) \qv (f\xi)(y\tto^h z,z\tto^k w) = \xi(x\tto^{hf}z,z\tto^k w) \]
\end{prop}
\epf
\begin{prop} \label{53prop5}
The morphism $\phi_P:P\to\,\up P$ which takes an object $a$ over $x$ to the component $[\la a , x \tto^\id x\ra]$ in $\comp(P/x)$
is a universal morphism to a dof. 
So the formula of Proposition~\ref{52prop2} gives indeed the reflection of $P$.
\end{prop}
\pf
First, we have to check that $\phi_P$ is indeed a morphism. By Proposition~\ref{53prop3},
for any $u:a\to b$ in $P$ over $f:x\to y$, 
\[ f[\la a , x \tto^\id x\ra] = [\la a , x \tto^f y\ra] = [\la b , y \tto^\id y\ra] \]
where the last equality follows from the arrow $u : \la b , y \tto^\id y\ra \to \la a , x \tto^f y\ra$ in $P/y$.
So the condition of Proposition~\ref{53prop1} is fulfilled. 

As for universality, given any morphism $\phi : P\to D$ to a dof, suppose that $\dbf\phi\circ\phi_p = \phi$.
Then $\phi a = \dbf\phi(\phi_p a) = \dbf\phi[\la a , x \tto^\id x\ra]$, so that
\[
\dbf\phi[\la a , x\tto^h y\ra] = \dbf\phi(h[\la a , x \tto^\id x\ra]) = h(\dbf\phi[\la a , x \tto^\id x\ra]) = h(\phi a) 
\]
which gives the unicity of $\dbf\phi$.
\[ 
\xymatrix@R=3.5pc@C=3.5pc{
P \ar[r]^{\phi_P}\ar[dr]_\phi  &  \up P \ar@{..>}[d]^{\dbf\phi}   \\ 
                               &  D                                }
\]
So we have to check that $\dbf\phi$ defined by
\[ \dbf\phi[\la a , x\tto^h y\ra] = h(\phi a) \]
is indeed a morphism $D\to P\rup$ in $\CatX$:
\begin{enumerate}
\item
$\dbf\phi$ is well-defined: if $u:a\to b$ in $P$ over $g:x\to z$ is an arrow $\la a , x\tto^h y\ra$ to $\la a , z\tto^k y\ra$
in $P/y$ (that is $kg = h$) then by Proposition~\ref{53prop1}
\[ h(\phi a) = kg(\phi a) = k(\phi b) \]
\item
$\dbf\phi$ is a morphism in $\CatX$ (or in $\dof X\simeq \Set^X$), since for any $f:y\to x$ in $X$
\[ \dbf\phi f[\la a , x\tto^h y\ra] = \dbf\phi [\la a , x\tto^{fh} z\ra] = fh(\phi a) = f\dbf\phi[\la a , x\tto^h y\ra] \]
\item
$\dbf\phi\circ\phi_p = \phi$, since for any $x\in Px$,
\(\, \dbf\phi(\phi_p a) = \dbf\phi([\la a , x \tto^\id x\ra]) = \id_x(\phi a) = \phi a \).
\end{enumerate}
\epf
\begin{prop} \label{53prop6}
The morphism $\phi_P:P\rup\,\to P$ which takes any $\xi\in(P\rup)x$ to the object $\xi(x\tto^\id x)$ in $P$ over $x$,
and any arrow $\la \xi , x\tto^f y \ra$ in $P\rup$ over $f$ to the arrow $\xi(x\tto^\id x,x\tto^f y)$ in $P$ over $f$
is a universal morphism from a dof. 
So the formula of Proposition~\ref{52prop1} gives indeed the coreflection of $P$.
\end{prop}
\pf
First, we have to check that $\phi_P$ is indeed a morphism in $\CatX$, that is that the conditions of Proposition~\ref{53prop2}
are fulfilled. But these follow directly from the corresponding ones about the $\xi$ in $(P\rup)x$ (see Proposition~\ref{53prop4}):
\begin{enumerate}
\item
\( \phi_P(\xi,x\tto^f y) = \xi(x\tto^\id x,x\tto^f y):\xi(x\tto^\id x)\to \xi(x\tto^f y)\,\); but \\ 
\( \xi(x\tto^\id x) = \phi_P\xi \,\,\) and \( \,\xi(x\tto^f y) = (f\xi)(y\tto^\id y) = \phi_P(f\xi) \,\) as required. 
\item
\( \phi_P(\xi,x\tto^\id x) = \xi(x\tto^\id x,x\tto^\id x) = \id_{\xi(\id_x)} = \id_{\phi_P\xi} \,\, \) and \\
\(
\phi_P(\xi,x\tto^{gf}z) = \xi(x\tto^\id x,x\tto^{gf}z) = \xi(x\tto^f y,y\tto^g z)\circ\xi(x\tto^\id x,f) \\
= (f\xi)(y\tto^\id y,y\tto^g z)\circ\xi(x\tto^\id x,f) = \phi_P(f\xi,y\tto^g z)\circ\phi_P(\xi,x\tto^f y)  
\)
\end{enumerate}
As for universality, given any morphism $\phi : D\to P$ from a dof, suppose that $\phi_P\circ\dbf\phi = \phi$.
Then \( \phi a = \phi_P(\dbf\phi a) = (\dbf\phi a)(x\tto^\id x) \), so that 
\[ (\dbf\phi a)(x\tto^h y) = (h(\dbf\phi a))(y\tto^\id y) = (\dbf\phi(ha))(y\tto^\id y) = \phi(ha) \]
and similarly $(\dbf\phi a)(x\tto^h y,y\tto^k z) = \phi(ha,y\tto^k z)$, giving the unicity of $\dbf\phi$.
\[ 
\xymatrix@R=3.5pc@C=3.5pc{
D \ar@{..>}[d]_{\dbf\phi}\ar[dr]^\phi  &   \\
P\rup \ar[r]^{\phi_P}                  &  P }
\]
So we have to check that $\dbf\phi$ defined by
\[ (\dbf\phi a)(x\tto^h y) = \phi(ha) \qv (\dbf\phi a)(x\tto^h y,y\tto^k z) = \phi(ha,y\tto^k z) \]
is indeed a morphism $D\to P\rup$ in $\CatX$:
\begin{enumerate}
\item
$\dbf\phi$ is well-defined, that is $\dbf\phi a$ is indeed a morphism $\up x \to P$, that is the conditions of Proposition~\ref{53prop4}:
are fulfilled. But these follow directly from the corresponding ones about $\phi$ (see Proposition~\ref{53prop2}):
\begin{enumerate}
\item
\( (\dbf\phi a)(x\tto^h y,y\tto^k z) = \phi(ha,y\tto^k z):\phi(ha)\to\phi(kha) \,\); but \\
\( \phi(ha) = (\dbf\phi a)(x\tto^h y) \,\,\) and \( \,\phi(kha) = (\dbf\phi a)(x\tto^{kh}z) \), as required.
\item
\( (\dbf\phi a)(x\tto^h y,y\tto^\id y) = \phi(ha,y\tto^\id y) = \id_{\phi(ha)} = \id_{(\dbf\phi a)h}\,\, \) and \\
\( 
(\dbf\phi a)(x\tto^h y,y\tto^{lk}w) = \phi(ha,y\tto^{lk}w) = \phi(kha,z\tto^l w)\circ\phi(ha,y\tto^k z) \\
= (\dbf\phi a)(x\tto^{kh}z,z\tto^l w)\circ(\dbf\phi a)(x\tto^h y,y\tto^k z) 
\)
\end{enumerate}
\item
$\dbf\phi$ is a morphism in $\CatX$ (or in $\dof X\simeq \Set^X$), since for any $f:x\to y$ in $X$ \\
\( \dbf\phi(fa)(x\tto^h y) = \phi(hfa) = (\dbf\phi a)(x\tto^{hf}y) = (f(\dbf\phi a))(x\tto^h y) \,\, \) and similarly \\ 
\(
\dbf\phi(fa)(x\tto^h y,y\tto^k z) = (f(\dbf\phi a))(x\tto^h y,y\tto^k z) 
\)
\item
$\phi_P\circ\dbf\phi = \phi$, since for any $a\in Dx$,
\( \,\phi_p(\dbf\phi a) = (\dbf\phi a)\id_x = \phi(\id_x a) = \phi a \,\),
and similarly for any arrow $\la a , x\tto^ f y\ra$ in $D$.
\end{enumerate}
\epf
\begin{remark} \label{53rmk1}
As suggested by Richard Wood, the reflection and coreflection of a category over $X$ in discrete fibrations can be 
factorized through fibrations over $X$.
Indeed, the inclusion $i:\dof X \inc \CatX$ factors (up to equivalence) through fibrations 
as follows
\[ \dof X \simeq \Set^X \inc \Cat^X \inc \CatX \]
where the last inclusion is given by the Grothendieck construction.
While the functor $(\disc)^X : \Set^X \inc \Cat^X$ clearly has both left and right adjoints 
\[ \comp^X\adj(\disc)^X\adj~\pt^X \] 
given \v pointwise" by components and points on $\Cat$, the same is true also for the Grothendieck functor.
Its left adjoint $\CatX\to\Cat^X$ is simply given by slicing: 
\[ P\mapsto P/-:X\to\Cat \]
(see also~\cite{her98}). 
Its right adjoint takes $P$ to the functor $\dbf P:X\to\Cat$ such that $\dbf P x$ is the category whose objects
are the morphisms $\xi : \,\up x \to P$ like in Proposition~\ref{53prop4}, while an arrow $\varphi : \xi \to \xi'$
is given by mappings $x\tto^h y\,\, \mapsto\, \varphi(x\tto^h y)$ such that
\begin{enumerate}
\item
$\varphi(x\tto^h y)\,$ is an arrow in $P$ over $\id_y$
\item
$\varphi(x\tto^h y):\xi(x\tto^h y)\to \xi'(x\tto^h y)$
\item
$\xi'(x\tto^h y,y\tto^k z)\circ\varphi(x\tto^h y) = \varphi(x\tto^{kh}z)\circ \xi(x\tto^h y,y\tto^k z)$.
\end{enumerate}
while for any $f:x\to x'$, the functor $\dbf P f:\dbf P x \to \dbf P x'$ is given by \\
\( (f\xi)(x'\tto^h y) = \xi(x\tto^{hf}y) \), \( (f\xi)(x'\tto^h y,y\tto^k z) = \xi(x\tto^{hf}y,y\tto^k z) \) on objects, and \\
\(\, (f\varphi)(x'\tto^h y) = \varphi(x\tto^{hf}y) \) on arrows.
\end{remark}

\section{Applications}
\label{applications}

In section~\ref{6sub1} we analyze reflections and coreflections in df's in relation with 
a \v change of base" given by a functor $X\to Y$, showing also their interdependence with Kan extensions.

The relevance of the reflection $\down(-)$ in df's over $X$ for the study of colimits of functors $\Phi:P\to X$ 
was shown thoroughly in~\cite{par73}. 
In section~\ref{6sub2}, we present those facts in a more conceptual way, resting on the universal property of $\down(-)$.
We begin with an unusual operative definition of (co)limits: the colimit of $\Phi\in\CatX$ is its reflection
in the full subcategory of \v principal" (or representable) df's: those with a terminal object.
Thus $\down(-):\CatX\to\df X$ can be seen as an intermediate step of this reflection.  

In section~\ref{4sub2} we have seen that any object $x\in X$ is an atom in $\CatX$.
In section~\ref{6sub3} we show that the same is true for any idempotent $e$ in $X$, 
and that the reflections of atoms are the retracts of representable functors.

In section~\ref{6sub4} we present \v another" context in which one can apply the reflection formula,
which allows in particular a direct and intuitive calculation of the \v best left approximation"
of a graph with an (idempotent, bijective or $n$-periodic) evolutive set.

Whenever opportune, we will be more careful about notations regarding categories over a base: 
referring to a category over $X$ we shall write e.g. $P\tto^p X$, rather than only $P$ as before.
Accordingly, its reflection will be denoted by $\down p$ rather than by $\down P$.

Dual statements, obtained by dualizing the base category, are often left to the reader.
Note that $\,\Phi\,\mapsto\,\Phi\op\,$ gives an isomorphism $\CatX\to\CatX\op$, such that
the following diagrams commute:
\[
\xymatrix@R=1.5pc@C=3pc{ 
{\df X} \ar[r]^-i \ar@2{-}[d]      & \CatX \ar@2{-}[d]      & {\dof X} \ar[l]_-j \ar@2{-}[d]  \\
{\dof{X\op}} \ar[r]^-j             & \CatX\op               & {\df{X\op}} \ar[l]_-i           }
\qq\qq
\xymatrix@R=1.5pc@C=3pc{ 
{\df X} \ar@2{-}[d]   & \CatX \ar@2{-}[d]\ar[l]_-{\down\,(-)}\ar[r]^-{\up\,(-)}   & {\dof X}  \ar@2{-}[d]  \\
{\dof{X\op}}          & \CatX\op  \ar[l]_-{\up\,(-)} \ar[r]^-{\down\,(-)}         & {\df{X\op}}            }
\]
Briefly, $\up\Phi\op = \,\down\Phi$ and $\down\Phi\op = \,\up\Phi$.
Furthermore, $\,\Lim\Phi\op = \Colim\Phi\,$, $\,\Colim\Phi\op = \Lim\Phi\,$,
$\,x/X\op = X/x\,$ and $\,X\op/x = x/X\,$.

\subsection{Kan extensions via reflection}
\label{6sub1}

Recall that any functor $X\to Y$ gives rise to a pair of adjoint functors 
\[ f_!\adj f^*:\CatY\to\CatX \] 
where $f_!$ is given by composition with $f$, while $f^*$ is obtained by the following pullback in $\Cat$:
\[ 
\xymatrix@R=3.5pc@C=3.5pc{
f^*P \ar[r] \ar[d]_{f^* p}  & P \ar[d]^ p   \\
X    \ar[r]^f                & Y               }
\]
Note that $\comp^X:\CatX\to\Set$ factors through $f_!:\CatX\to\CatY$ and $\comp^Y$ up to isomorphisms;
indeed, the components of a category over a base are given by those of the total category (see section~\ref{3sub3}) and $P$ and $f_!P$ 
have the same total category. As a consequence also $\disc_X$ factors through $f^*$ and $\disc_Y$ (up to isomorphisms):
\eq  \label{61eq2}
\xymatrix@R=3.5pc@C=3.5pc{
\CatX   \ar[r]^{f_!} \ar[d]^\comp    & \CatY \ar[d]^\comp   \\
\Set    \ar@2{-}[r]                  & \Set }
\qq\qq
\xymatrix@R=3.5pc@C=3.5pc{
\CatX                          & \CatY \ar[l]_-{f^*}   \\
\Set \ar[u]_\disc \ar@2{-}[r]  & \Set  \ar[u]_\disc }
\eeq
Furthermore, the pair $f_!\adj f^*$ satisfies the Frobenius law, that is the morphism
\eq  \label{61eq1}
f_!\pi_1\wedge(\eps_Q\circ f_!\pi_2) : f_!(P\times f^*Q) \to f_! P\times Q 
\eeq
is an isomorphism for any $P\in\CatX$ and $Q\in\CatY$, $\eps$ being the counit of the adjunction. 
This follows essentially from the fact that the composition of the two pullback squares below, 
expressed by $f_!(P\times f^*Q)$, gives another pullback, expressed by $f_! P\times Q$: 
\[
\xymatrix@R=3.5pc@C=3.5pc{
P\times f^*Q\ar[r]\ar[d] & f^*Q \ar[r] \ar[d]_{f^*q} & Q\ar[d]^q   \\
P\ar[r]^p                &  X    \ar[r]^f            & Y               }
\]
As it is easily checked (see~\cite{law70}), this implies that $f^*$ preserves existing exponentials.
In particular, by~(\ref{61eq2}), $f^*$ preserves complements:
\eq   \label{61eq0}
\xymatrix@R=3.5pc@C=3.5pc{
\CatX    \ar[d]^{(\disc S)^-}    & \CatY \ar[l]_{f^*}\ar[d]^{(\disc S)^-}  \\
(\CatX)\op                     &  (\CatY)\op\ar[l]_{f^*}                } 
\eeq 
The following adjunction-like property holds:
\begin{prop}  \label{61prop1}
For any functor $f:X\to Y$ there is a bijection
\[ \ten_X(P,f^*Q) \cong \ten_Y(f_!P,Q) \]
natural in $P\in\CatX$ and $Q\in\CatY$.
\end{prop}
\pf 
We have 
\eqa*   
\ten_X(P,f^*Q) = \comp^X(P\times f^*Q) \cong \comp^Y(f_!(P\times f^*Q)) \\
\cong \comp^Y(f_!P\times Q) = \ten_Y(f_!P,Q) 
\eeqa*
where we have used the first of~(\ref{61eq2}) and the Frobenius law.
\epf
\begin{corol}  \label{61prop12}
Let $f:X\to Y$ any functor. Then for any atom $T\tto^t X$ in $\CatX$, $f_! t:T\to Y$ is an atom in $\CatY$.
\end{corol}
\pf
We have the following chain of natural bijections:
\[ \ten_{\CatY}(f_!x,D) \cong \ten_{\CatX}(x,f^*D) \cong \hom_{\CatX}(x,f^*D) \cong \hom_{\CatY}(f_!x,D) \]
\epf
Now we turn to Kan extensions of set functors. 
First recall that pullbacks preserve df's and dof's: 
if $D$ over $Y$ is in $\dof Y$ then $f^*D$ is in $\dof X$.
So $f^*:\CatY\to\CatX$ restricts to $f^*:\df Y\to\df X$ and $f^*:\dof Y\to\dof X$:
\[
\xymatrix@R=3.5pc@C=3.5pc{
{\df X} \ar[d]^i & {\df Y}\ar[l]_{f^*}\ar[d]^i \\ 
\CatX            & \CatY\ar[l]_{f^*}                 }
\qq\qq
\xymatrix@R=3.5pc@C=3.5pc{
{\dof X} \ar[d]^j & {\dof Y}\ar[l]_{f^*}\ar[d]^j \\ 
\CatX             & \CatY\ar[l]_{f^*}                 }
\]
In terms of presheaves $f^*:\Set^Y\to\Set^X$ can be defined by composition:
\[ f^*D = D\circ f : X\to\Set \]
\begin{remark}   \label{61rmk1}
Thus also diagram~(\ref{61eq0}) restricts to 
\eq   \label{61eq01}
\xymatrix@R=3.5pc@C=3.5pc{
\Set^{X\op}    \ar[d]^{\neg_S} & \Set^{Y\op} \ar[l]_{f^*}\ar[d]^{\neg_S}  \\
(\Set^X)\op                     &  (\Set^Y)\op\ar[l]_{f^*}                } 
\eeq 
E.g., if $f:X\to Y$ is a quotient obtained by adding \v laws" (commutative diagrams) to $X$,
then $f^*$ is full and faithful and we can say that complementation in $\CatX$ \v preserves (or better, reverses) algebraic laws".
For example, let $X$ be the free category on $\xymatrix@1{x \ar@<-1ex>[r]\ar@<1ex>[r] & y \ar[l]}$, 
and let $D\in\dof X$ be a dof which happens to be a reflexive graph. 
Then its complements (in particular $(\neg D)(1+1)$ given by $\xymatrix@1{\P(Dy) \ar@<-1ex>[r]\ar@<1ex>[r] & \P(Dx) \ar[l]}$)
are cylinders (mappings with two sections); conversely, the complements of cylinders are reflexive graphs (see~\cite{law89}).
\end{remark}
We are interested in the left and right adjoints of $f^*$:
\[ \exists_f \adj f^* \adj \forall_f :\dof X \to \dof Y \]
also known as Kan extensions. We begin with the following
\begin{lemma}  \label{61prop2}
For any functor $f:X\to Y$ and any $P$ over $X$,
\[ \exists_f(\up P) \cong \,\up(f_!P) \]
\end{lemma}
\pf
Since in the diagram below the square of right adjoint arrows commutes up to isomorphisms, 
so does also the square of left adjoint arrows: 
\eq
\xymatrix@R=4pc@C=4pc{
{\CatX} \ar@<1ex>[r]^{f_!}_{\bot}\ar@<-1ex>[d]_{\up\,(-)}^{\adj} & {\CatY}\ar@<1ex>[l]^{f^*}\ar@<-1ex>[d]_{\up\,(-)}^{\adj} \\ 
{\dof X} \ar@<1ex>[r]^{\exists_f}_{\bot}\ar@<-1ex>[u]_{j}        & {\dof Y}\ar@<1ex>[l]^{f^*}\ar@<-1ex>[u]_{j} }
\eeq
\epf
\begin{corol}  \label{61prop3}
The following formulas hold: 
\[ \exists_f D = \,\up(f_!D) \qv \up f = \exists_f \id_X \]
\end{corol}
\pf
For any $D\in\dof X$ and $f:X\to Y$,
\( \,\exists_f D = \exists_f \up D = \,\up(f_!D) \) \\
while for any category $X\tto^f Y$ over $Y$,
\[ \up f = \,\up(f_!\id_X) = \exists_f(\up\id_X) = \exists_f \id_X \]
since $\id_X:X\to X$ is the terminal df $\,\disc 1$. 
\epf
\begin{corol}  \label{61prop4}
The left Kan extension of a set functor $D:X\to\Set$ along $f:X\to Y$ is given by the \v coend formula"
\[ (\exists_f D)y = Y(f-,y)\otimes D = \Coend\,(Y(f-,y)\times D-) \]
\end{corol}
\pf
By the above corollary, $\exists_f D = \,\up(f_!D)$ so that
\[ (\exists_f D)y = \,\up(f_!D)y = \ten_{\CatY}(\down y,f_! D) = \ten_{\CatX}(f^*\down y,D) = \Coend\,(Y(f-,y)\times D-) \]
where Proposition~\ref{61prop1} and property~(\ref{43eq2}) have been used.
\epf
\begin{remark}
Supposing on the contrary that the coend formula of Corollary~\ref{61prop4} is known, 
one gets again the reflection formula using Corollary~\ref{61prop3} as follows:
\[
\begin{array}{c}
(\up f)y  \\ \hline
(\exists_f \id_X)y  \\ \hline
\ten_X(f^*\down y,\id_X)  \\ \hline
\comp((f^*\down y)\times\id_X)  \\ \hline
\comp(f^*\down y)  \\ \hline
\comp(f\times\down y)  \\ \hline
\comp(f/y)
\end{array}
\]
\end{remark}
\begin{remark}
As it was the case for the {\em coreflection} formula of Proposition~\ref{52prop1},
the \v end formula" for {\em right} Kan extensions can be obtained directly as follows
(we use Remark~\ref{35rmk1} in the last step):  
\[
\begin{array}{c}
(\forall_f D)y \\ \hline
\Nat(\up y,\forall_f D) \\ \hline 
\Nat(f^*\up y,D) \\ \hline 
\End\,\Set(Y(y,f-),D- ) 
\end{array}
\]
\end{remark}

\subsection{Limits and colimits}
\label{6sub2}

Given a functor $\Phi:P\to X$ and an object $x\in X$, a cone from $\Phi$ to $x$
is a family of arrows $\, f_i:\Phi i \to x\, , i\in P\,$, such that $\Phi u\circ f_i = f_j$ for any $u:i\to j$ in $P$.
If $f:x\to y$ is an arrow in $X$, the family $\, f\circ f_i:\Phi i \to y\,$ a cone with vertex $y$.
So we have a functor
\[ \Cone(\Phi,-) : X\to\Set \]  
This functor can be substantiated in two ways, that lead to two different points of view on limits and colimits.
The usual way is to see a cone as a natural transformation from $\Phi$ to $\Delta x$, so that one defines
\[ \Cone(\Phi,x) := X^P(\Phi,\Delta x) \]
But a cone from $\Phi$ to $x$ is also the object mapping $\, i\,\mapsto\, f_i \,$
of a functor $\Phi\to X/x$ over $X$ (see Proposition~\ref{52prop1}), 
and we can define the above functor equivalently as 
\[ \Cone(\Phi,x) := \CatX(\Phi,X/x) \]
Thus, while for any fixed $\Phi$ we get isomorphic functors
\eq  \label{62eq0}
\CatX(\Phi,X/-) \cong \Nat(\Phi,\Delta -) 
\eeq
we in fact have two different modules \v cone"
\[ \Cone:X^P\to X \qv \Cone:\CatX\to X \]
both representable on the right.

Since a colimit of $\Phi$ is a representing object for $\Cone(\Phi,-)$, we have also
two different operative definitions of the colimit {\em functor}:
either as the partially defined left adjoint 
\[ \Colim : X^P\to X \]
to the functor $\,\Delta:X\to X^P\,$, or as the partially defined left adjoint 
\[ \Colim : \CatX\to X \]
to the functor $\, X/- : X \to \CatX\,$.
Dually, we have two limit functors 
\[ \Lim : X^P\to X \qv \Lim : \CatX\to X\op \]
respectively right adjoint to $\,\Delta:X\to X^P\,$ and left adjoint to $\, -/X : X\op\to \CatX\,$.
\begin{remarks}
\begin{itemize}
\item
Since $X/-$, unlike $\Delta$, is always full and faithful, the counit $\eps_x:\Colim X/x\to x$
of the adjunction $\Colim\adj X/-$ is an isomorphism. 
(Indeed, this follows also from the fact that $X/x$ has a terminal object over $x$, that is $x\tto^\id x$.)
The unit $\eta_\Phi:\Phi\to X/\Colim\Phi$ is again the universal cone.
\item
The two colimit functors (which agree on each single $\Phi:P\to X$ as an object of two different categories)
can be merged in a single
\[ \Colim : \CatX^*\to X \]
(see~\cite{mac71},~\cite{par73} and~\cite{kan58}) where the \v supercomma category" $\CatX^*$ is the Gro\-then\-dieck category
associated to the 2-representable functor $\Cat(-,X):\Cat\to\Cat$ (which includes both $\CatX$ and all the $X^P$ 
as subcategories). Indeed, we have the further module \v cone" (which includes the other ones)
\[ \Cone(\Phi,x) := \CatX^*(\Phi,\delta x) \] 
where $\delta:X\to\CatX^*$ is the full and faithful functor which sends an object $x\in X$ 
to the corresponding category over $X$. 
\item
If $H=\hom_X(F-,-)$, then $H$ over $X$ is the forgetful functor $U$ from $F$-algebras to $X$ 
(see the last of Examples~\ref{32ex1}); then the strong dinatural transformations $H\to X/x$
are the cones $U\to x$, as in~\cite{din}.  
\end{itemize}
\end{remarks}
Since $X/x$ is a df (associated to the representable $X(-,x)$),
the functor $X/-$ factors through $i:\df X\inc\CatX$. Thus its left adjoint $\Colim$
also factors through $\down(-):\CatX\to\df X$.
More explicitly, we have the following chain of isomorphisms:
\[ \Cone(\Phi,-) \cong \CatX(\Phi,X/-) \cong \CatX(\down\Phi,X/-) \cong \Cone(\down\Phi,-) \]
so that $\Colim\,\Phi$ exists iff $\Colim\down\Phi$ exists, and if this is the case
\[ \Colim\,\Phi \,\cong\, \Colim\down\Phi \]
where $\down\Phi:P'\to X$ is the reflection of $\Phi:P\to X$ in discrete fibrations;
dually \(\, \Lim\,\Phi \,\cong\, \Lim\up\Phi \,\), either one existing if the other one does. 
\begin{prop}  \label{62prop1}
The colimit of a functor depends only on its reflection in df's.
Dually, the limit of a functor depends only on its reflection in dof's.
\end{prop}
\epf
The above considerations can be parametrized:
given a functor $f:X\to Y$, one can compose the adjunctions $f_!\adj f^*$ and $\Colim\adj Y\-$ as follows:
\[
\xymatrix@C=3.5pc{ \CatX \ar@<1ex>[r]^{f_!}_\bot\ar@/^2pc/@<1ex>[rr]_{\Colim(f\circ-)}  
& \CatY \ar@<1ex>[l]^{f^*} \ar@<1ex>[r]^\Colim_\bot  &  Y \ar@<1ex>[l]^{Y/-}\ar@/^2pc/@<1ex>[ll]_{f/-} } 
\]
So we have the partially defined adjunction with parameter $f:X\to Y$:
\eq  \label{62eq1}
\Colim(f\circ-) \adj f/-: Y\to\CatX 
\eeq
(Alternatively, one easily checks directly that $\,\Cone(f\circ\Phi,y) \cong \CatX(\Phi,f/y)\,$.)
Since the category $\, f/y = f^*(Y/y)\,$ over $X$ is the df associated to $Y(f-,y):X\op\to\Set$,
the functor $f/-$ factors through $i:\df X\inc\CatX$. 
Thus its left adjoint $\Colim(f\circ-)$ also factors through $\down(-):\CatX\to\df X$.
More explicitly, we have
\eqa* 
\Cone(f\circ\Phi,-) \cong \CatX(f\circ\Phi,Y/-) \cong \CatX(\Phi,f/-) \cong \CatX(\down\Phi,f/-) \\
\cong \CatX(f\circ\down\Phi,Y/-) \cong \Cone(f\circ\down\Phi,-) 
\eeqa*
Thus $\Colim\,f\circ\Phi$ exists iff $\Colim\,f\circ\down\Phi$ exists, and if this is the case
\[ \Colim(f\circ\Phi) \cong \Colim(f\circ\down\Phi) \]
Furthermore, the weak colimits of $f\circ\Phi$ coincide with those of $f\circ\down\Phi$.
In particular we have proved in a very direct manner the following proposition: 
\begin{prop}  \label{62prop2}
If two functors $\Phi:P\to X$ and $\Psi:Q\to X$ with the same codomain have isomorphic 
reflections $\,\down\Phi\,\, \cong\,\, \down\Psi\,$, then for any $f:X\to Y$
\[ \Colim(f\circ\Phi) \,\cong\, \Colim(f\circ\Psi) \]
either side existing if the other one does.
Furthermore, $f\circ\Phi:P\to Y$ and $f\circ\Psi:Q\to Y$ have the same weak colimits.
\end{prop}
\epf
To prove the converse of Proposition~\ref{62prop2}, we begin with the following
\begin{lemma}  \label{62prop21}
Let $P\tto^\Phi X$ be a category over $X$, $D$ a dof over $X$, and $\y$ 
the \v Yoneda embedding" $\,\, -/X:X\op\to\dof X\,$. Then there is a bijection 
\[ \CatX(\Phi,D) \,\cong\, \Cone(\y\circ\Phi\op,D) \]
natural in $D\in\dof X$.
\end{lemma}
\pf
Rephrasing Proposition~\ref{53prop1}, a morphism $\phi:\Phi\to D$ in $\CatX$ is given by
a family of morphisms $\phi_i:\Phi i\to D,\, i\in P$, in $\CatX$, such that for any $u:i\to j$ in $P$,
$\Phi u:\phi_i(\Phi i)\,\mapsto\,\phi_j(\Phi j)$ in $D$.
By the universality of $\up(\Phi i) = \Phi i/X$, the latter corresponds to a family of morphisms
$\phi_i:\Phi i/X\to D,\, i\in P$, in $\CatX$, such that for any $u:i\to j$ in $P$,
$\Phi u:\phi_i(\Phi i\tto^\id\Phi i)\,\mapsto\,\phi_j(\Phi j\tto^\id\Phi j)$ in $D$;
that is $\phi_i(\Phi i\tto^{\Phi u}\Phi j) = \phi_j(\Phi j\tto^\id\Phi j)$, as required by the cone
condition, since $(\y\circ\Phi\op)u:\Phi j/X\to\Phi i/X:\Phi j\tto^\id\Phi j\mapsto\Phi i\tto^{\Phi u}\Phi j$.
\epf
\begin{corollary}  \label{62prop3}
The cone $\,\y\circ\Phi\op\to\,\up\Phi\,$ in $\dof X$, corresponding to the unit $\eta:\Phi\to\,\up\Phi$ 
of the reflection $\up(-):\CatX\to\dof X$ in dof's, is a universal cone. 
Briefly, $\up\Phi$ is the colimit of $\,\y\circ\Phi\op:P\op\to\dof X$ 
(or the limit of $\,\y\op\circ\Phi:P\to(\dof X)\op$).
\end{corollary}
\epf
Putting together Proposition~\ref{62prop2} and Corollary~\ref{62prop3}, we get (see~\cite{par73}):
\begin{corollary}  \label{62prop4}
Let $\Phi:P\to X$ and $\Psi:Q\to X$ be two functors with the same codomain.
The following are equivalent:
\begin{itemize}
\item
\( \up\Phi\, \cong\,\, \up\Psi \,\);
\item
\(\, \Lim(f\circ\Phi) \cong \Lim(f\circ\Psi) \,\),
for any $f:X\to Y$, either side existing if the other one does.
\end{itemize}
\end{corollary}
\epf
\begin{remarks}   \label{62rmk1}
\begin{itemize}
\item
If in Lemma~\ref{62prop21} we take $\,D = x/X = \y x$, then we have a bijection  
\[ \Cone(\Phi\op,x) \,\cong\, \Cone(x,\Phi) \,\cong\, \CatX(\Phi,x/X) \,\cong\, \Cone(\y\circ\Phi\op,\y x) \]
which is the same that we would obtain by applying directly the Yoneda embedding $\,\y : X\op\to\dof X\,$
to the cones $\,\Phi\op \to x$.
\item
Alternatively, following~\cite{par73}, one can prove Corollary~\ref{62prop3} by noting that
since colimits in $\Set^{X\op}$ are computed pointwise, 
\eq  \label{62eq3}
(\Colim(\y\circ\Phi\op))x \cong \Colim((\y\circ\Phi\op)x) \cong \Colim\, X(\Phi-,x) \cong \comp(\Phi/x) \cong (\up\Phi)x 
\eeq
Conversely, by Corollary~\ref{62prop3}, the~(\ref{62eq3}) above gives a new proof of the reflection formula.  
\item
A further proof of Corollary~\ref{62prop3} can be obtained noticing that morphisms $\Phi\to D$ in $\CatX$
correspond also to cones $1\to D\circ\Phi$ in $\Set$, considering $D$ as a presheaf.
(This fact also allows another proof of the second assertion in Proposition~\ref{61prop3},
as explained in ~\cite{law70}.)
Then we have
\[
\begin{array}{c}
\dof X(\Colim(\y\circ\Phi\op),D) \\ \hline
\Lim\,\dof X(\y\circ\Phi\op -,D) \\ \hline
\Lim(D\circ\Phi) \\ \hline
\Cone(1,D\circ\Phi) \\ \hline
\CatX(\Phi,D)
\end{array}
\]
\end{itemize}
\end{remarks}
Among the consequences of Corollary~\ref{62prop4}, the following are particularly significant (see~\cite{par73}):
\begin{corol}  \label{62prop5}
\begin{enumerate}
\item
A category over a base and its reflections $\,\up\Phi$ and $\,\down\Phi$ have the same components.
\item
The functor $\Phi$ is initial iff $\,\,\up\Phi \,\cong\, \id_X\,$.
\item
The unit $\, \eta:\Phi\to\up\Phi \,$ is an initial functor between the respective total categories.
\item
The object $x\in X$ is the absolute limit of $\Phi$ iff $\,\,\up\Phi\,\cong\, x/X\,$.
\item
The object $x\in X$ is an absolute weak limit of $\Phi$ iff $\,\up\Phi$ 
is a retract of $\,x/X\,$ in $\dof X$.
\end{enumerate}
\end{corol}
\pf
\begin{enumerate}
\item
By Remark~\ref{34rmk1}, for any $P\tto^p X$ in $\CatX$ we have
\[ \comp P \cong \Colim \Delta 1 \cong \Colim (\Delta 1\circ p) \cong \Colim (\Delta 1\circ\down p) \cong \comp(\down P)  \] 
The same result follows from the fact that since in the diagram below the square of right adjoint arrows commutes up to isomorphisms, 
so does also the square of left adjoint arrows:
\[
\xymatrix@R=4pc@C=4pc{
{\CatX} \ar@<1ex>[r]^{\comp}_{\bot}\ar@<-1ex>[d]_{\down\,(-)}^{\adj} & {\Set}\ar@<1ex>[l]^{\disc}\ar@2{-}[d] \\ 
{\df X} \ar@<1ex>[r]^{\comp}_{\bot}\ar@<-1ex>[u]_{i}           & {\Set}\ar@<1ex>[l]^{\disc}             }
\]
\item
Apply Corollary~\ref{62prop4} with $\Psi = X\tto^\id X$, the terminal df.
\item
By point 2 just proved, we must show that for any dof $D$ over $\up P$ and any $\phi:\eta\to d$ over $\up P$,
$\phi$ factors uniquely over $\up P$ as in the diagram below:
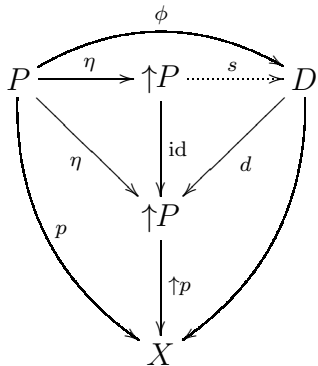
\[
\xymatrix@R=3pc@C=3pc{ 
P \ar[r]^\eta \ar@/^1.5pc/[rr]^\phi \ar[dr]_\eta \ar@/_1.5pc/[ddr]^p &  \up P \ar[d]^\id \ar@{..>}[r]^s  & D \ar[dl]^d \ar@/^1.5pc/[ddl]                 \\
                                                                     &  \up P \ar[d]^{\up\,\,p}          &                    \\
                                                                     &  X                                &                     }
\] 
Since $\,\up p\circ d\,$ is a dof over $X$, and by the universality of $\eta$ over $X$, 
$\,\phi:p\to\,\,\up p\circ d\,$ factors over $X$ as $\,\phi = s\circ\eta$. 
Then $d\circ s\circ\eta \cong d\circ\phi \cong \eta$, so that $d\circ s \cong \id_{\up\,P}$, again
by the universality of $\eta$ over $X$. 
Thus $s$ is a morphism over $\up\Phi$, giving the existence of the desired factorization.
As for unicity, another such $s$ would give another factorization through $\eta$ of $\phi$ over $X$.
(Note that this proof is not specifically linked to the situation now considered, 
unlike those given in~\cite{par73},~\cite{her98}, and~\cite{kel82}.)  
\item
Apply Corollary~\ref{62prop4} with $\Psi = 1\tto^x X$.
Note that each cone $f_i:x\to\Phi i$ in $X$, corresponding to $\Phi\to x/X$ in $\CatX$, induces by universality
\eq  \label{62eq4}
\xymatrix@R=3pc@C=3pc{                         & \up\Phi \ar@{..>}[d]^\alpha \\
                      \Phi\ar[r] \ar[ur]^\eta  & x/X                          }
\eeq
Recalling Corollary~\ref{62prop3} and the first of Remarks~\ref{62rmk1}, we may see $\alpha$ as the universally induced morphism 
\[ \Colim(\y\circ\Phi\op) \to x/X \]
in $\CatX$, through which the cone $\y f_i = f_i/X : \Phi i/X \to x/X$ factors.
Then the cone $f_i:x\to\Phi i$ in $X$ is an absolute limit iff the cone $f_i/X :\Phi i/X \to x/X$ is a limit in $\dof X$
(or in $\CatX$) iff $\alpha$ is an isomorphism.
\item
If the cone $f_i:x\to\Phi i$ in $X$ is an absolute weak limit in $X$, then the cone $f_i/X:\Phi i/X \to x/X$ is a weak limit.
Thus the universally induced $\alpha$, as in~(\ref{62eq4}) above, has a retraction. 
The reverse implication will be proved in~\ref{63rmk3}. 
\end{enumerate}
\epf
\begin{remarks}
\begin{itemize}
\item
If $Y=\Set$, the adjunction~(\ref{62eq1}) becomes the adjunction 
\[ \ten(A,-) \adj \neg A : \Set\to\CatX \]
of section~\ref{4sub4}, with parameter $A$ in $\df X$.
\item
Since $f/y$ is a df over $X$, the adjunction~(\ref{62eq1}) restricts to
\eq  \label{62eq2}
\Colim(f\circ-) \adj f/-: Y\to\df X\simeq\Set^{X\op} 
\eeq
which displays the weighted colimit construction (as a functor of the weight $A\in\Set^{X\op}$)
\[ A \,\mapsto\, A*f \in Y \]
as the left adjoint to the functor
\[ y \,\mapsto\, Y(f-,y) \in \Set^{X\op} \]
(of course the $A$ in $\Colim(f\circ A)$ is intended over $X$ so that the colimit is taken
on the \v category of elements" of the corresponding presheaf).
Thus, at least in the non-enriched context, the paradigmatical adjunction~(\ref{62eq2})
(as illustrated e.g. in~\cite{mam91}) is nothing but a particular case of the more general~(\ref{62eq1}).  
(See also~\cite{mac71}, where the restriction of~(\ref{62eq0}) to df's is referred to as the \v coyoneda lemma" on page~62.)

Note that from Proposition~\ref{62prop2} we have
\[ \Colim(f\circ\Phi) \cong \Colim(f\circ\down\Phi) = \,\down\Phi*f \] 
\end{itemize}
\end{remarks}

\subsection{Atoms in $\CatX$ and their reflections}
\label{6sub3}

We begin by recalling some well-known facts about idempotents and their splitting in $\Set$ and $\Set^X$. 
Let $\e$ be the monoid which represents the idempotent arrows of categories, 
that is the one whose unique non-identity arrow is idempotent.
\begin{prop}  \label{63prop0}
For any category $X$ and any functor $\,\e\tto^e X$, the following are equivalent:
\begin{enumerate}
\item
The idempotent $e:x\to x$ in $X$ splits.
\item 
The functor $e$ has a limit.
\item 
The functor $e$ has a colimit.
\item 
The pair $e,\id:x\to x$ has an equalizer.
\item 
The pair $e,\id:x\to x$ has a coequalizer.
\end{enumerate}
If this is the case, all the above limits and colimits are absolute and canonically isomorphic.
\end{prop}
\epf
Now let $e:D\to D$ an idempotent mapping, and $\,\e\tto^e \Set$ the corresponding functor.
As a dof over $\e$, the category $e$ is (apart from the identities) exactly what we would call the 
\v graph of the endomapping $e$".
\begin{corol}  \label{63prop01}
For any functor $\,\e\tto^e \Set$, its limit, given by the $\fix\, e$ 
(the set of fixed points of the endomapping $e:D\to D$),
and its colimit, given by the components of $e$ over $\e$, are canonically isomorphic. 
\end{corol}
\epf
\noindent Explicitly, the bijection $\beta :\Lim\,e \to \Colim\,e$ is given by $\beta x = [x]$ and $\beta\inv [x] = ex$. 
\begin{corol}  \label{63prop1}
Any idempotent $\,\e\tto^e X$ in $X$ is an atom in $\CatX$.
\end{corol}
\pf
For any dof $D$, $\hom(\e,D)$ gives $\fix De$, while $\ten(\e,D)$ 
gives the components $\e\times D$ which, over $\e$, is the graph of $De:Dx\to Dx$ as in Corollary~\ref{63prop01}.
Dually, $e$ is also a left atom.
\epf
\begin{prop}  \label{63prop11}
For any atom $\, T\tto^t X$ of $\CatX$, 
\[ \down T \otimes - \,\cong\, \Nat(\up T,-):\Set^X\to\Set \qv -\,\otimes\up T \,\cong\, \Nat(\down T,-):\Set^{X\op}\to\Set \] 
For any idempotent atom $\,\e\tto^e X$,  
\[ \down e \otimes D \cong \Nat(\up e, D) \cong \fix De \qv  A\,\otimes\up e \cong \Nat(\down e, A) \cong \fix Ae \]
\end{prop}
\pf
Since $T$ is a right atom, for any dof $D$ one has the following bijections, natural in $D\in\dof X$:
\[ \q \Nat(\up T,D) \cong \hom(\up T, D) \cong \hom(T, D) \cong \ten(T, D) \cong \ten(\down T, D) \cong\, \down T \otimes D \]
and similarly for any df $A$.
The second part follows from Corollary~\ref{63prop1}.
\epf
\begin{remarks}  \label{63rmk1}
\begin{enumerate}
\item
In particular, for any idempotent $\, e:x_0\to x_0$ in $X$,  
\[ \hom(\up e, \up x) \cong \ten(\down e, \up x) \cong \fix X(x,e) \]
where $X(x,e) = e\circ- : X(x,x_0)\to X(x,x_0)$.
But the reflection formula (or co-Yoneda) gives $\down e \cong \ten(\down e, \up -)$, so that
\[ \down e \cong \hom(\up e, \up -) \cong \fix X(-,e) \]
that is, $\down e$ is the subfunctor of $\,\down x_0 = X(-,x_0)\,$ given by the arrows fixed by composition with $e$:
\[ x\tto^f x_0 \in \fix X(x,e) \iff e\circ f = f \]
Dually,
\[ \up e \cong \hom(\down e, \down -) \cong \fix X(e,-) \]
\item
In particular, for any idempotents $\, e:x\to x$ in $X$ and $\, e':y\to y$,  
\[ \hom(\up e, \up e') \cong \ten(\down e, \up e') \cong \fix (\up e')e \]
is the subset of $(\up e')x$ fixed by $(\up e')e$. 
By the above remark, this is the set of arrows in $\fix X(e',x)$ fixed also by $e$,
that is the arrows $f:y\to x$ such that $e\circ f = f = f\circ e'$. 
Dually \(\, \hom(\up e, \up e')\, \) is the set of arrows $f:x\to y$ such that $f\circ e = f = e'\circ f$.
\end{enumerate}
\end{remarks}
\noindent By Yoneda, any idempotent in $\Set^X$ on a representable $X(x_0,-)$ has the form 
\[ X(e,-):X(x_0,-)\to X(x_0,-) \]
for a unique idempotent $\, e:x_0\to x_0$ in $X$.
\begin{corol}  \label{63prop02}
Let $D:X\to\Set$. The following are equivalent:
\begin{enumerate}
\item
$D$ is a retract in $\Set^X$ of the representable functor $X(x_0,-)$, associated to the idempotent $X(e,-):X(x_0,-)\to X(x_0,-)$.
\item
$D$ is isomorphic to the subfunctor $\fix X(e,-)$ of $X(x_0,-)$.
\item
$D$ is isomorphic to the quotient of $X(x_0,-)$ given by \(\, f\circ e \sim f \).
\item
$D$ is the reflection of the idempotent atom $e$: $\, D\,\cong\,\,\up e$.
\end{enumerate}
\end{corol}
\pf
By Proposition~\ref{63prop0}, the retract $D$ is the limit of $\e\ttto^{X(e,-)} \Set^X$; 
but the latter is computed pointwise, so that 1 is equivalent to 2 and 3 by Corollary~\ref{63prop01}.
That 4 is equivalent to 2 follows from the first of Remarks~\ref{63rmk1}; 
alternatively note that the functor in 3 is $\comp (\e/-)$, and the latter is $\up e$ by the reflection formula.
\epf
\begin{prop}  \label{63prop12}
Let $A$ be a df and $D$ a dof on $X$ such that
\[ A\otimes - \cong \Nat(D,-):\Set^X\to\Set \]
Then $D$ is a retract in $\Set^X$ of a representable functor.
\end{prop}
\pf
Let $u\in A\otimes D = \comp(A\times D)$ be an universal element of $A\otimes -$. Then $u$ is the component 
of a pair $\la a,d \ra$ over, say, $x_0\in X$: 
\[ u = [\la a,d \ra],\q a\in Ax_0,\, d\in Dx_0 \]
Let $\iota:D\to X(x_0,-)$ be the unique morphism in $\Set^X$ such that $(A\otimes\iota)u =[\la a,\id_{x_0} \ra]$:
\[ [\la a,\id_{x_0} \ra] = (A\otimes\iota)u = \comp(A\times\iota)[\la a,d \ra] = [\la a,\iota d \ra] \] 
and let $\rho:X(x_0,-)\to D$ be the unique morphism in $\Set^X$ such that $\rho\,\id_{x_0} = d$. Then 
\[ A\otimes\rho\circ\iota = \comp(A\times\rho\circ\iota): [\la a,d \ra] \,\mapsto\, [\la a,d \ra] \]
so that $\rho\circ\iota = \id_D$, as required.
Note that $D$ splits the idempotent $\iota\circ\rho:X(x_0,-)\to X(x_0,-)$ which corresponds to
the idempotent $e = \iota d:x_0\to x_0$ in $X$. 
\epf
\begin{prop}  \label{63prop2}
Let $A$ be a df and $D$ a dof on $X$. The following are equivalent:
\begin{enumerate}
\item
$A\otimes - \cong \Nat(D,-):\Set^X\to\Set$.
\item
$-\otimes D \cong \Nat(A,-):\Set^{X\op}\to\Set$.
\item
$A\,\cong\,\,\down e$ and $D\,\cong\,\,\up e$, for an idempotent atom $e$ in $\CatX$.
\item
$A\,\cong\,\,\down T$ and $D\,\cong\,\,\up T$, for an atom $T$ in $\CatX$.
\end{enumerate}
\end{prop}
\pf
Trivially 3 implies 4, and 4 implies 1 and 2 by Proposition~\ref{63prop11}.   
From Proposition~\ref{63prop12} and Corollary~\ref{63prop02} it follows that 1 implies that $D\,\cong\,\up e$.
Furthermore, 
\[ A\,\cong\, A\,\otimes \up - \,\cong\,\Nat(D,\up -)\,\cong\,\Nat(\up e,\up -)\,\cong\,\,\down e\,\otimes\up - \,\cong\,\,\down e \]
where Proposition~\ref{63prop11} have been used again. Thus 1 implies 3. 
By duality we also have the equivalence between 2 and 3 (note that 3 is autodual).
\epf
\begin{remark}
In~\cite{kel05} it is shown, in an enriched context, the equivalence between 1 and 2 and the fact that $A$ and $D$ are adjoint modules.
\end{remark}
\begin{corol}
Any right atom is also a left atom, and vice versa.
\end{corol}
\pf
Suppose that $T$ is a right atom.
As in Proposition~\ref{63prop11}, \(\, \down T \otimes - \,\cong\, \Nat(\up T,-) \,\),
and so by Proposition~\ref{63prop2} \(\, -\,\otimes\up T \,\cong\, \Nat(\down T,-) \,\). 
Then, for any df $A$,
\[ \ten(A,T) \,\cong\, A\,\otimes\up T \,\cong\, \Nat(\down T,A) \,\cong\, \hom(T,A) \]
that is, $T$ is also a left atom.
\epf
\begin{corol}  \label{63prop3}
Let $A$ be a df on $X$. The following are equivalent:
\begin{enumerate}
\item
$A$ is a retract of a representable functor.
\item
$A \,\cong\,\,\down e$ for an idempotent atom $e$ in $\CatX$.
\item
$A \,\cong\,\,\down T$ for an atom $T$ in $\CatX$. 
\item
$A\otimes -:\Set^X\to\Set$ is representable.
\item
$\Nat(A,-)\cong -\otimes D:\Set^{X\op}\to\Set$ for a $D\in\Set^X$.
\item
$A\otimes -:\Set^X\to\Set$ has a left adjoint.
\item
$\Nat(A,-):\Set^{X\op}\to\Set$ has a right adjoint.
\item
$A\otimes - :\Set^X\to\Set$ preserves all limits.
\item
$\Nat(A,-):\Set^{X\op}\to\Set$ preserves all colimits.
\item
$A$ is a left adjoint module.
\end{enumerate}
\end{corol}
\pf
We have already seen the equivalence between the first five properties.
If $A\otimes - \cong \Nat(D,-)$, then $(-)\cdot D\adj A\otimes -$.
If $\Nat(A,-)\cong -\otimes D$, then $\Nat(A,-)\adj\neg D$.
For the equivalence between the last three properties and 4 and 5, we refer again to~\cite{kel05},
where they are proved in an enriched context. 
\epf
\begin{remark}  \label{63rmk3}
We can now complete the proof of Proposition~\ref{62prop5}, showing that if $\up\Phi$ is a retract of a
representable functor $X(x,-)$ then $x$ is an absolute weak limit of $\Phi$.
By Corollary~\ref{63prop02}, $\up\Phi \cong\, \up e$, for an idempotent $e:x\to x$.
But $e:x\to x$ itself is a cone from $x$ to $e:\e\to X$, and it is weakly universal:
\[
\xymatrix@R=2pc@C=2pc{ y \ar@{..>}[d]_f \ar[dr]^f  &                    \\
                       x \ar[r]^e                  &  x \ar@(ur,dr)^e    }
\]
Furthermore, this weak limit is clearly absolute. 
Then, by Proposition~\ref{62prop2}, $x$ is also an absolute weak limit of $\Phi$. 
\end{remark}
By Proposition~\ref{63prop2}, for any atom $T\tto^t X$ there exists an idempotent $e:x\to x$ in $X$ 
such that $\up T = \,\up e$ and $\down T = \,\down e$.
This idempotent is not unique, since e.g. if $e$ splits as $x\tto^r y\tto^i x$, then $\up e = \,\up y$.
Anyway, by the above remark, any atom has an absolute weak (co)limit.
\begin{corol}
For an atom $T\tto^t X$ the following are equivalent:
\begin{itemize}
\item
$T$ has a limit.
\item
$T$ has a colimit.
\item
$T$ has an absolute limit.
\item
$T$ has an absolute colimit.
\end{itemize}
\end{corol}
\pf
Indeed, by Proposition~\ref{63prop0}, all the above are equivalent for any $e:\e\to X$, and so 
by Corollary~\ref{62prop4} they are equivalent also for $T$.
\epf
\begin{remark}
If the conditions of the above proposition are satisfied, we could say that $t$ \v converges" to its (co)limit.
In particular, any split idempotent in $X$ converges to its retracts, and a category is Cauchy complete
iff any atom converges. Furthermore, any functor $f:X\to Y$ is \v continuous": if $t$ converges to $x$,
then $f_! t$ converges to $fx$ (see Corollary~\ref{61prop12}). 
\end{remark}

\noindent Recall that the Cauchy completion of a category can be obtained as the full subcategory of $\Set^{X\op}$ 
generated by the retracts of representable functors (see e.g.~\cite{bor94}).
But we have seen that the latter have the form $\down T$ or $\down e$, and so we get the first part of the
following proposition.
The second part follows from the second of Remarks~\ref{63rmk1}.
\begin{prop}  
The reflections of atoms (or of idempotent atoms) in df's generate the Cauchy completion of $X$.
Furthermore, given two idempotent arrows $e:x\to x$ and $e':y\to y$ in $X$, $\hom(\down e,\down e')$ 
is the set of arrows $f:x\to y$ such that $f\circ e = f = e'\circ f$.
So the Cauchy completion is the same thing as the \v Karoubi envelope" of $X$ 
{\rm (see~\cite{rey04}, \cite{lam86} and~\cite{law89}).} 
\end{prop}
\epf

\subsection{Graphs and evolutive sets}
\label{6sub4}

Let $X$ be a graph in $\Grph$. A graph $P$ over $X$ is a (graph) opfibration if for any $f$ in $X$,
the fiber $Pf$ (defined like in section~\ref{3sub1}) \v is" a mapping.
Fibrations are defined dually. So we have the inclusions
\[ i : \df X \inc \GrphX \qv j : \dof X \inc \GrphX \]
of fibrations and opfibrations in graphs over $X$. 
These have left and right adjoints
\[ \down(-) \adj i \adj (-)\rdown \qv \up(-) \adj j \adj (-)\rup \]
as can be shown in two ways:
\begin{itemize}
\item
Since $\Grph$ is a presheaf category, $\GrphX$ is also a presheaf category, namely on the total category of the df $X$; 
the above inclusions may be seen as induced by functors from this category, so that the reflection and coreflection 
are given by the left and the right Kan extensions along these functors.  
\item
The functor $\,\GrphX \to \CatFX\,$ given by
\[ P\ttto^p X \q\mapsto\q {\cal F}\!P\ttto^{{\cal F}\! p}{\cal F}\!X \]
where $\,{\cal F}:\Grph\to\Cat\,$ is the \v free category" functor,
is a full and faithful functor whose values are (up to isomorphisms) the \v UFL functors" to ${\cal F}\!X$ (see e.g.~\cite{bun00});
and $p$ is a graph fibration iff ${\cal F}\! p$ is a df. 
So the reflection of the category ${\cal F}\!P\ttto^{{\cal F}\! p}{\cal F}\!X$ over ${\cal F}\!X$ in df's 
(being itself a UFL functor) gives also the reflection of $P\ttto^p X$ in graphs fibrations:
\eq
\xymatrix@R=4pc@C=4pc{
{\GrphX} \ar[r]\ar[d]_{\down\,(-)} & {\CatFX}\ar[d]^{\down\,(-)} \\ 
{\df X} \ar[r]^\sim                & {\df{{\cal F}X}}               }
\eeq
\end{itemize}
But perhaps it is more interesting to see why and how the reflection formulas can be applied directly to this context:
\begin{itemize}
\item
Every node $x\in X$, as a graph over $X$, is a strong atom in $\GrphX$ since
\[ \hom(x,P) \cong \ten(x,P) \cong Px \]
(note that the fiber $Px$ is {\em always} discrete).
\item
The reflections $\up x$ and $\down x$ are given by $(\up x)y = {\cal F}\!X(x,y)$ and $(\down x)y = {\cal F}\!X(y,x)$.
\item
Any graph over $X$ has a complement (since $\GrphX$ is cartesian closed); like in section~\ref{3sub1}, one sees that
the complement of a fibration is valued in opfibrations, and vice versa.
\item
The above facts are all what is needed to prove the reflection formulas
\[ (\down P)x = \ten(\up x,P) \qv (\up P)x = \ten(\down x,P) \]
as was done for categories over a base.  
Of course, one has also the coreflection formulas
\[ (P\rdown)x = \hom(\down x,P) \qv (P\rup)x = \hom(\up x,P) \]
\end{itemize}
In particular, if $X$ is the terminal graph (the \v loop"), we have the inclusions
\[ i : \Lendo \inc \Grph \qv j : \Rendo \inc \Grph \]
of anti-evolutive sets and evolutive sets (or left and right endomappings) in graphs.
The reflections of the dot atom $\d$ (see section~\ref{4sub2}) are the infinite chain and anti-chain: 
\[ 
\up \d = \q \xymatrix@1@C=2pc{
{\bullet}\ar[r] & {\bullet}\ar[r] & {\bullet}\ar[r] & {\bullet}\ar[r] & {\bullet}\ar@{..>}[r] &  }
\]
\[ 
\down \d = \q \xymatrix@1@C=2pc{
{\bullet} & {\bullet}\ar[l] & {\bullet}\ar[l] & {\bullet}\ar[l] & {\bullet}\ar[l] & \ar@{..>}[l] }
\]
So the reflection $\up P$ of a graph $P$ in evolutive sets is given by the set $\comp(\down \d \times P)$,
with the action given by the right shift of the anti-chain $\down\d$:
\[ [n,x]\,\mapsto\, [n+1,x] \]
\begin{examples}  \label{64ex1}
\begin{enumerate}
\item
Let $P$ be the graph
\( \q \xymatrix@1@C=2pc{{\bullet}\ar[r] & {\bullet} }\q \)
which is not a (right) endomapping because there are no arrows out of one node.
Applying the reflection formula, we multiply $P$ with the anti-chain $\down \d$, getting the graph
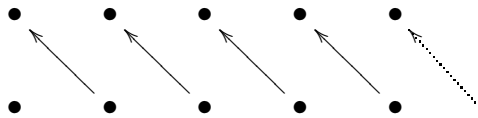
\[
\xymatrix@R=2pc@C=2pc{
{\bullet} & {\bullet}        & {\bullet}        & {\bullet}        & {\bullet}        &                 \\
{\bullet} & {\bullet}\ar[ul] & {\bullet}\ar[ul] & {\bullet}\ar[ul] & {\bullet}\ar[ul] &  \ar@{..>}[ul]  } 
\]
with an infinite number of components, which are the nodes of the reflection of $P$.
Furthermore, the action on it is given by translation, which sends any component in the one on its right.
So $\up P$ is the chain \( \, \xymatrix@1@C=2pc{
{\bullet}\ar[r] & {\bullet}\ar[r] & {\bullet}\ar[r] & {\bullet}\ar[r] & {\bullet}\ar@{..>}[r] &  } \, \),
wherein the missing codomains in $P$ have been added.
On the other hand, there are no chains in $P$ (that is, $\hom(\up\d,P) = \emptyset$), so that the coreflection $P\rup$ 
is the void endomapping: the nodes with no codomains have been deleted.
\item
Let $P$ be the graph
\( \qq \xymatrix@1@C=2pc{{\bullet}\ar@(ul,dl) & {\bullet}\ar[l]\ar[r] & {\bullet}\ar@(ur,dr) }\qq \)
which is not a (right) endomapping because there are two arrows out of one node.
Applying the reflection formula, we multiply $P$ with the antichain $\down \d$, getting the two-components graph
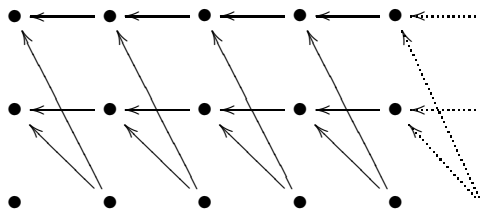
\[
\xymatrix@R=2pc@C=2pc{
{\bullet} & {\bullet}\ar[l] & {\bullet}\ar[l] & {\bullet}\ar[l] & {\bullet}\ar[l] & \ar@{..>}[l]    \\
{\bullet} & {\bullet}\ar[l] & {\bullet}\ar[l] & {\bullet}\ar[l] & {\bullet}\ar[l] & \ar@{..>}[l]    \\
{\bullet} & {\bullet}\ar[ul]\ar[uul] & {\bullet}\ar[ul]\ar[uul] & {\bullet}\ar[ul]\ar[uul] & {\bullet}\ar[ul]\ar[uul] & \ar@{..>}[ul]\ar@{..>}[uul]  } 
\]
So the reflection of $P$ has two nodes. 
Furthermore, the action on it is given again by right translation, so that the second component is a fixed point and we get
\[ 
\up P = \xymatrix@1@C=2pc{{\bullet}\ar[r] & {\bullet}\ar@(ur,dr)  }
\]
wherein the multiple codomains in $P$ have been identified.
On the other hand, there are four chains in $P$ so that $P\rup$ has four nodes:
\[ 
P\rup\q = \qq\xymatrix@1@C=2pc{{\bullet}\ar@(ul,dl) & {\bullet}\ar[l] & {\bullet}\ar[r] & {\bullet}\ar@(ur,dr) }
\]
that is, the nodes with multiple codomains have been splitted.
\item
If \( P = \, \xymatrix@1@C=2pc{{\bullet}\ar@/^1pc/[rr] & {\bullet}\ar[l]\ar[r] & {\bullet} }\, \), $\down\d\times P$ is the graph 
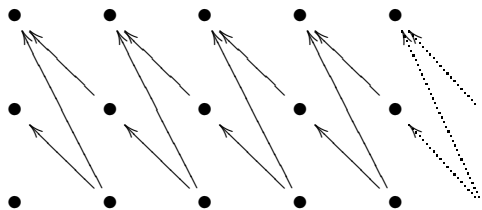
\[
\xymatrix@R=2pc@C=2pc{
{\bullet} & {\bullet}        & {\bullet}        & {\bullet}        & {\bullet}        &               \\
{\bullet} & {\bullet}\ar[ul] & {\bullet}\ar[ul] & {\bullet}\ar[ul] & {\bullet}\ar[ul] & \ar@{..>}[ul]    \\
{\bullet} & {\bullet}\ar[ul]\ar[uul] & {\bullet}\ar[ul]\ar[uul] & {\bullet}\ar[ul]\ar[uul] & {\bullet}\ar[ul]\ar[uul] & \ar@{..>}[ul]\ar@{..>}[uul]  } 
\]
so that \(\, \up P = \,\xymatrix@1@C=2pc{{\bullet}\ar[r] & {\bullet}\ar@(ur,dr)  }\q\,\, \) again.
One should compare this example with the technique used in~\cite{rey04} to compute the same reflection.
\item
An example where both the phenomena of adding and identifying codomains, are present is given by the graph 
\( \q \xymatrix@1@C=2pc{{\bullet} & {\bullet}\ar[l]\ar[r] & {\bullet} }\q \), whose reflection is again the chain.
\item
If \( \, P_1 = \q\,\, \xymatrix@1@C=2pc{{\bullet}\ar@(ul,dl)\ar[r] & {\bullet} } \,\),
\( \, P_2 = \q\,\, \xymatrix@1@C=2pc{{\bullet}\ar@(ul,dl)\ar[r] & {\bullet}\ar@(ur,dr) }\q\,\, \)
and \( \, P_3 = \q\,\, \xymatrix@1@C=2pc{{\bullet}\ar@(ul,dl)\ar@(ur,dr) }\q\,\, \),
then $\up P_1$, $\up P_2$ and $\up P_3$ are all the loop (identification prevails in $\up P_1$). 
As for coreflections, $P_1\rup$ is the loop, $P_2\rup$ is the sum of a loop and a left infinite chain ending with a loop on the right,
and $P_3\rup$ is the set $2^\N$ of sequences in a two-element set, under the action $f\,\mapsto\, f(1+ -)$. 
\end{enumerate}
\end{examples} 
$\Rendo$ and $\Lendo$ can be defined as the full subcategories of $\Grph$ generated by the objects orthogonal 
respectively to 
\eq   \label{64eq1}
\xymatrix@1@C=2pc{
\n\, =\, {\circ}\ar[r] & {\bullet}\ar[r] & {\bullet}\ar[r] & {\bullet}\ar@{..>}[r] &  } 
\qv
\xymatrix@1@C=2pc{
\n\op\, =\, {\circ} & {\bullet}\ar[l] & {\bullet}\ar[l] & {\bullet}\ar[l] & \ar@{..>}[l] }
\eeq
meaning the morphism of graphs from the dot $\d$ to the chain (respectively, the anti-chain) 
which sends the unique node of $\d$ to the highlighted node. 
Similarly, orthogonality with respect to
\eq  \label{64eq3}
\ee\, =\, \xymatrix@1@C=2pc{{\circ}\ar[r] & {\bullet}\ar@(ur,dr) }  \qq\qv \ee\op\, = \, \xymatrix@1@C=2pc{{\circ} & {\bullet}\ar[l]\ar@(ur,dr) }
\eeq
\eq  \label{64eq4}
\xymatrix@1@C=2pc{
\z\, =\q \ar@{..>}[r] & {\bullet}\ar[r] & {\bullet}\ar[r] & {\circ}\ar[r] & {\bullet}\ar[r] & {\bullet}\ar@{..>}[r] &  }
\eeq
\eq  \label{64eq5}
\z_n \, =\q \xymatrix@R=1.5pc@C=1.5pc{  \circ \ar[r] & \bullet \ar@/^/@{..>}[dl] \\  
                                      \bullet \ar[u] &                            } \q (n\,\, {\rm arrows})
\eeq
defines respectively the full subcategories $\Lidem$ and $\Ridem$ of left and right idempotent endomappings, 
and the full subcategories $\Bij$ of bijective endomappings and $\Bij_n$ of n-periodic bijective endomappings.
So we have the following inclusions:
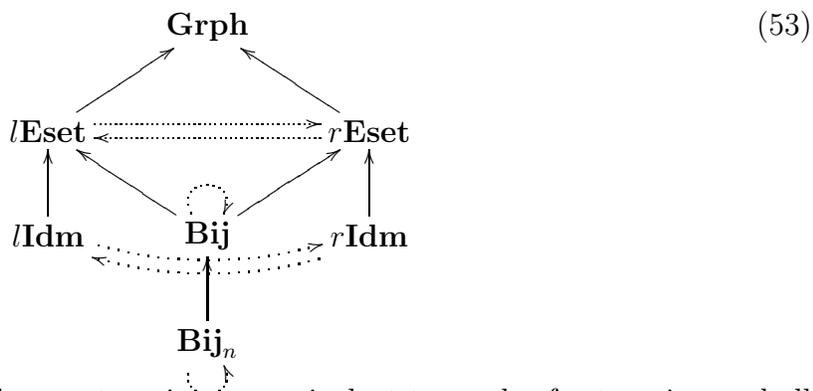
\eq  \label{64eq6}
\xymatrix@R=2pc@C=2pc{                            & \Grph                                 &                                    \\
            \Lendo\ar[ur] \ar@<.5ex>@{..>}[rr]    &                                       & \Rendo\ar[ul] \ar@<.5ex>@{..>}[ll] \\
            \Lidem\ar[u]  \ar@<.5ex>@/_1pc/@{..>}[rr]    &  \Bij\ar[ul]\ar[ur] \ar@{..>}@(ul,ur) & \Ridem\ar[u] \ar@<.5ex>@/^1pc/@{..>}[ll]  \\
                                                  &  \Bij_n\ar[u]\ar@{..>}@(dl,dr)        &                                      }
\eeq
\begin{remark}  \label{64rmk1}
All the above categories are equivalent to presheaf categories, and all the above inclusions correspond to 
functors induced by epimorphism between the corresponding categories of shapes.
So they have left and right adjoints given by Kan extensions.
But also in this case, a direct application of the formulas gives easier computations of these adjoints. 
\end{remark}
\begin{itemize}
\item
The effect of complementation in $\Grph$ is displayed by the dotted arrows in diagram~(\ref{64eq6}):
the complement of a left (right) idempotent endomapping is valued in right (left) idempotent endomappings;
the complement of a bijective endomapping is valued in bijective endomappings;
the complement of a $n$-periodic endomapping is valued in $n$-periodic endomappings.
These facts may be seen as a consequence of Remark~\ref{61rmk1} and Remark~\ref{64rmk1}.
\item
Since the~(\ref{64eq1}), (\ref{64eq3}), (\ref{64eq4}) and~(\ref{64eq5}) are all orthogonal to themselves,
they serve also as the reflections of the dot graph $\d$ in the corresponding subcategories
(where the highlighted morphism from $\d$ now gives universal arrow). 
\end{itemize}
\noindent The above facts are all what is needed to prove the reflection formulas of graphs in each of the above subcategories;
e.g., \(\, \ten(\ee\op,P) = \comp(\ee\op \times P) \,\) gives the nodes of the reflection of the graph $P$ in (right) idempotent mappings.
Of course, one has also the corresponding coreflection formulas;
e.g., \(\, \Grph(\ee,P) \,\) gives the nodes of the coreflection of the graph $P$ in (right) idempotent mappings.
\begin{examples}
\begin{itemize}
\item
It is easy to see that if $P$ is connected (that is, $\comp P = 1$) then $\ee\op\times P$ has $n+1$ components,
where $n$ is the number of nodes that are not codomains of any arrow. So the reflection in (right) idempotent mappings
acts on each component mantaining such nodes, and collapsing the rest of it to the fixed point. 
\item
The reflection $\up P = \,\down P$ of a graph $P$ in bijective endomappings is obtained by taking the components of $P\times\z$.
For the graphs of Examples~\ref{64ex1} we obtain either $\z$ itself or the fixed point $\z_1$.
If \(\, P\, =\, \xymatrix@R=2pc@C=2pc{ \bullet \ar[r] & \bullet \ar@/^/[r] & \bullet \ar@/^/[l] } \,\), 
then $\up P$ and $P\rup$ are both $\z_2$. 
\item
The reflection $\up P = \,\down P$ of a graph $P$ in $n$-periodic endomappings is obtained by taking the components of $P\times\z_n$.
Since, as is easily checked, 
\[ \z_k\times \z_n = \gcd(k,n)\cdot \z_{\lcm(k,n)} \] 
we deduce that $\up \z_k$ has $\gcd(k,n)$ nodes and so, being connected by Corollary~\ref{62prop5}, 
\[ \up \z_k = \z_{\gcd(k,n)} \] 
Any bijective endomapping $P$ is a sum of cycles: $P = \sum_{k=1}^\infty S_k\cdot \z_k$, so that 
\[ \up P =  \up\Big(\sum_{k=1}^\infty S_k\cdot \z_k \Big) = \sum_{k=1}^\infty S_k\cdot \up \z_k = \sum_{k=1}^\infty S_k\cdot \z_{\gcd(k,n)} \]
\end{itemize}
\end{examples} 
\begin{remarks}
\begin{itemize}
\item
The graphs~(\ref{64eq3}) and~(\ref{64eq5}) are instances of the following
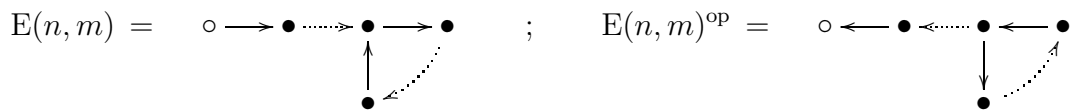
\[
\ee(n,m) \, =\q \xymatrix@R=1.5pc@C=1.5pc{ \circ\ar[r] & \bullet\ar@{..>}[r] & \bullet\ar[r]   & \bullet\ar@{..>}@/^/[dl] \\  
                                                         &                   & \bullet\ar[u]   &                           } 
\qv                                                      
\ee(n,m)\op \, =\q \xymatrix@R=1.5pc@C=1.5pc{ \circ & \bullet\ar[l] & \bullet\ar@{..>}[l]\ar[d]   & \bullet\ar[l]      \\  
                                                    &               &  \bullet\ar@{..>}@/_/[ur]   &                     }                        
\]
where $n\geq 1$ is the length of the cycle, while $m\geq 0$ is the length of the chain.
As before, these graphs allow to compute the reflection (and coreflection) in the corresponding full subcategories of $\Grph$,
with the usual formulas.
\item
Obviously, when the graph $P$ grows bigger, it becomes more and more difficult to visualize the components of $\down\d\times P$.
So an appropriate software would be helpful. 
Note that if $P$ is finite, then also when $\down\d$ is infinite, only a finite portion of the product is needed to capture the reflection.
The size of this portion depends on that of $P$ itself.   
\end{itemize}
\end{remarks}

\begin{refs}

\bibitem[Borceux, 1994]{bor94} F. Borceux (1994), {\em Handbook of Categorical Algebra 1 (Basic Category Theory)}, 
Encyclopedia of Mathematics and its applications, vol. 50, Cambridge University Press.

\bibitem[Bunge \& Niefield, 2000]{bun00} M. Bunge and S. Niefeld (2000), Exponentiability and Single Universes, 
{\em J. Pure Appl. Algebra} {\bf 148}, 217-250.
 
\bibitem[Ghani, Uustalu \& Vene, 2004]{din} N. Ghani, T. Uustalu and V. Vene (2004), {\em Build, Augment and Destroy. Universally},
Proceeding of the Asian Symposium on Programming Languages, LNCS 3302, Springer, Berlin, 327-347.

\bibitem[Hermida, 1998]{her98} C. Hermida (1998), A Note on Comprehensive Factorization, expository note, 
available at {\tt http://maggie.cs.queensu.ca/chermida/ }

\bibitem[Higgins, 1971]{hig71} P.J. Higgins (1971), {\em Notes on categories and grupoids}, 
Van Nostrand Reinhold Company.

\bibitem[Kan, 1958]{kan58} D.M. Kan (1958), Adjoint Functors, {\em Trans. Am. Math. Soc.} {\bf 87}, 294-329.

\bibitem[Kelly, 1982]{kel82} G.M. Kelly (1982), {\em Basic Concepts of Enriched Category Theory}, 
London Mathematical Society Lecture Note Series, 64, Cambridge University Press.
Republished in {\em Reprints in Theory and Appl. Cat.} No. 10 (2005) 1-136.
 
\bibitem[Kelly \& Schmitt, 2005]{kel05} G.M. Kelly and V. Schmitt (2005), Notes on Enriched Categories with Colimits of some Class, 
{\em Theory and Appl. Cat.} {\bf 17}, 399-423.

\bibitem[Lambek \& Scott, 1986]{lam86} J. Lambek and P.J. Scott (1986), {\em Introduction to Higher Order Categorical Logic}, 
Cambridge studies in advanced mathematics, vol. 7, Cambridge University Press.

\bibitem[La Palme, Reyes \& Zolfaghari, 2004]{rey04} M. La Palme Reyes, G.E. Reyes, H. Zolfaghari (2004), 
{\em Generic Figures and their Glueings}, Polimetrica S.a.s., Monza. 

\bibitem[Lawvere, 1970]{law70} F.W. Lawvere (1970), {\em Equality in Hyperdoctrines and the Comprehension Scheme as an Adjoint Functor},
Proceedings of the AMS Symposium on Pure Mathematics, XVII, 1-14. 

\bibitem[Lawvere, 1973]{law73} F.W. Lawvere (1973), Metric Spaces, Generalized Logic and Closed Categories,
{\em Rend. Sem. Mat. Fis. Milano} {\bf 43}, 135-166. Republished in {\em Reprints in Theory and Appl. Cat.} No. 1 (2002) 1-37.

\bibitem[Lawvere, 1986]{law86} F.W. Lawvere (1986), Taking Categories seriously, 
{\em Revista Colombiana de Matematicas} {\bf 20}, 147-178. Republished in {\em Reprints in Theory and Appl. Cat.} No. 8 (2005) 1-24.

\bibitem[Lawvere, 1989]{law89} F.W. Lawvere (1989), {\em Qualitative Distinctions between some Toposes of Generalized Graphs},
Proceedings of the AMS Symposium on Categories in Computer Science and Logic, Contemporary Mathematics, vol. 92, 261-299.

\bibitem[Lawvere, 1996]{law96} F.W. Lawvere (1996), {\em Adjoints in and among Bicategories},
Proceedings of the 1994 Conference in Memory of Roberto Magari,
Logic \& Algebra, Lectures Notes in Pure and Applied Algebra, 180:181-189,
Ed. Ursini Aglian\`o, Marcel Dekker, Inc. Basel, New York. 

\bibitem[Mac Lane, 1971]{mac71} S. Mac Lane (1971), {\em Categories for the Working Mathematician}, 
Graduate Texts in Mathematics, vol. 5, Springer, Berlin.

\bibitem[Mac Lane \& Moerdijk, 1991]{mam91} S. Mac Lane and I. Moerdijk (1991), {\em Sheaves in Geometry and Logic (a First Introduction to Topos Theory)}, 
Universitext, Springer, Berlin.


\bibitem[Par\'e, 1973]{par73} R. Par\'e (1973), Connected Components and Colimits,
{\em J. Pure Appl. Algebra} {\bf 3}, 21-42.

\bibitem[Par\'e \& Roman, 1998]{par98} R. Par\'e and L. Rom\'an (1998), Dinatural Numbers,
{\em J. Pure Appl. Algebra} {\bf 128}(1), 33-92.

\bibitem[Pisani, 2005]{bip} C. Pisani (2005), Bipolar Spaces, preprint, math.CT/0512194

\end{refs}

\end{document}